\documentclass[11pt]{article}
\usepackage{graphicx} % Required for inserting images
\graphicspath{{media/}}
\usepackage{amsmath}
\usepackage{amssymb}

\graphicspath{{images/}}

\newcommand{\bfb}{{\bm b}}

\newcommand{\bfx}{{\bm x}}
\newcommand{\bfy}{{\bm y}}
\newcommand{\bfz}{{\bm z}}
\newcommand{\bfu}{{\bm u}}

\newcommand{\bfr}{{\bm r}}

\newcommand{\bfv}{{\bm v}}

\newcommand{\bfzero}{{\bm 0}}

\def\true{_{_{\mathrm{true}}}}
\def\inv{_{_{\mathrm{inv}}}}
\def\opt{_{_{\mathrm{opt}}}}
\def\gcv{_{_{\mathrm{gcv}}}}
\def\wgcv{_{_{\mathrm{wgcv}}}}
\def\disc{_{_{\mathrm{disc}}}}

\newcommand{\bfeps}{{\bm \varepsilon}}

\newcommand{\CD}{{\cal D}}

\newcommand{\CO}{{\cal O}}

\newcommand{\CG}{{\cal G}}

\newcommand{\R}{\ensuremath{\mathds{R}}}

\makeatletter
\newcommand*\bigldot{\mathpalette\bigldot@{.5}}
\newcommand*\bigldot@[2]{\mathbin{{\hbox{\scalebox{#2}{$\m@th#1\bullet$}}}}}
\makeatother

\makeatletter
\newcommand*\bigdec{\mathpalette\bigdec@{.5}}
\newcommand*\bigdec@[2]{\mathbin{{\hbox{\scalebox{#2}{$\!\!\m@th#1\bullet$}}}}}
\makeatother

\makeatletter
\newcommand*\bigcdot{\mathpalette\bigcdot@{.5}}
\newcommand*\bigcdot@[2]{\mathbin{\vcenter{\hbox{\scalebox{#2}{$\m@th#1\bullet$}}}}}
\makeatother
\def\dotmult{\,\bigldot *\;}

\def\ttdotmult{$\,\bigldot *\;$}

\usepackage{theorem,ifthen,algorithm,algorithmic}
\usepackage{amssymb,amsfonts,latexsym,dsfont}
\usepackage{tikz,pgflibraryplotmarks}
\usepackage{fullpage}
\usepackage{graphicx}
\usepackage{subcaption}
\usepackage{mathrsfs}
\usepackage{algorithm}
\usepackage{algorithmic}
\usepackage{amsmath}   
\usepackage{bm}
\usepackage{float}
\usepackage{caption}
\usepackage{hyperref}
\usepackage{epstopdf}
\usepackage{comment}
\usepackage{setspace}

\title{Half-Precision Kronecker Product SVD Preconditioner \\ for Structured Inverse Problems}

\author{
Yizhou Chen\thanks{Department of Mathematics, Emory University, Atlanta, GA 30322, USA. E-mail: rileychen111@gmail.com} \and
Xiang Ji\thanks{Department of Mathematics, Emory University, Atlanta, GA 30322, USA. Corresponding author. \newline E-mail: zoejix@outlook.com} \and
James Nagy\thanks{Department of Mathematics, Emory University, Atlanta, GA 30322, USA. E-mail: jnagy@emory.edu. Research supported in part by the U.S. National Science Foundation, grant codes DMS-2038118 and DMS-2208294.}}

\begin{document}

\maketitle

\begin{abstract}
In this paper we investigate the use of half-precision Kronecker product singular value decomposition (SVD) approximations as preconditioners for large-scale Tikhonov regularized least squares problems. 
Half precision reduces storage requirements and has the potential to greatly speedup computations on certain GPU architectures.  We consider both standard PCG and flexible PCG algorithms, and investigate, through numerical experiments on image deblurring problems, the trade-offs between potentially faster convergence with the additional cost per iteration when using this preconditioning approach. Moreover, we also investigate the use of several regularization parameter choice methods, including generalized cross validation and the discrepancy principle. 
\end{abstract}

\section{Introduction}
Inverse problems can occur in various fields, when observations are acquired to infer the original data. They can be described as the following linear system:

\begin{equation}\label{invp}
    \bfb = A\bfx\true + \bfeps
\end{equation}
where $\bfb$ is the outside measurements containing some noise $\bfeps$ and $A$, usually ill-conditioned and large-scale, represents some mapping of the true solution $\bfx\true$ \cite{afkham2021learning}. Since the noise is unknown, we can only calculate an approximate solution $\bfx$, instead of $\bfx\true$. A na\"ive approach computes an approximation by ignoring the noise, and either solving $A\bfx=\bfb$, if $A$ is $n\times n$ and non-singular, or solving the least squares problem, 
\begin{equation}\label{lsp}
    \min_\bfx \|\bfb-A\bfx\|^2_2\,.
\end{equation}

In either case, when $A$ is large, iterative methods are needed. An additional difficulty is that $A$ is typically very ill-conditioned, with singular values clustering at $0$. In this case, noise in the measured data is greatly amplified when computing the na\"ive solution, resulting in very poor approximations of $\bfx\true$. To get around this difficulty, regularization methods are usually incorporated into the solution process.

One of the most common approaches is Tikhonov regularization, which is implemented to balance the smoothness of the solution and the residual term. It can be written as the least squares problem, 
\begin{equation}\label{lsp:Tik}
   \min_\bfx \|\bfb-A\bfx\|^2_2+\lambda^2\|\bfx\|^2_2
    =\min_\bfx \left\|\left[ \begin{array}{c}
\bfb\\
\bfzero\\
\end{array} \right]-\left[ \begin{array}{c}
A\\
\lambda I\\
\end{array} \right]\bfx\right\|^2_2 \, . 
\end{equation}
The regularization parameter, $\lambda$, balances the residual norm with the solution norm \cite{golub1999tikhonov}; large $\lambda$ favors small norm (smoother) solutions, while small $\lambda$ favors large residual norms (and hence, less smooth solutions). There are various ways to choose the regularization parameter $\lambda$, such as generalized cross validation (GCV), weighted GCV, and the discrepancy principle \cite{chung2008weighted,hansen2010discrete}.

The conjugate gradient (CG) algorithm is a well-known and efficient iterative method for solving large-scale sparse (or structured) linear systems. 
For least squares problems, CG can be applied to the normal equations, which for the Tikhonov regularization problem (\ref{lsp:Tik})takes the form
\begin{equation*}
    \left[ \begin{array}{cc}
A^T&\lambda I\\
\end{array} \right]\left[ \begin{array}{c}
A\\
\lambda I\\
\end{array} \right]\bfx = \left[ \begin{array}{cc}
A^T&\lambda I\\
\end{array} \right]\left[ \begin{array}{c}
\bfb\\
\bfzero\\
\end{array} \right] \, .
\end{equation*}
\begin{equation}
\label{tik:normal1}
(A^TA+\lambda^2I)\bfx=A^T\bfb\, .
\end{equation}
Normally, it is recommended to use either CGLS or LSQR implementations \cite{bjorck1990least}. The overall cost of this algorithm depends on the number of iterations and the cost of matrix-vector multiplications involving $A$ and $A^T$, which can be done efficiently if $A$ is sparse or if it has an exploitable structure (e.g., Toeplitz or Kronecker product).

Preconditioning can be used in the CG algorithm (PCG) to accelerate convergence. The basic idea, for equation (\ref{tik:normal1}), is to construct a matrix $M$ that has three properties. The first is that $M \approx A^TA+\lambda^2 I$, such that the eigenvalues of $M^{-1}(A^TA+\lambda^2I)$ are clustered around a point bounded away from zero. The second is that solving the linear system $M\bfz=\bfr$ is not too expensive, which is required at each iteration of the PCG algorithm. The third property is that it is not too expensive to construct the approximating matrix $M$. 

There exists various ways to construct $M$. In this paper we use Kronecker product singular value decomposition (SVD) approximations. 
To further accelerate the process, we use mixed-precision computations in the PCG algorithm. Specifically, we use 16 bit half-precision arithmetic to store $M$ and to compute solutions of $M\bfz=\bfr$, and 64 bit double precision for the rest of the computations.
Half precision reduces storage requirements and has the potential to greatly speedup computations on certain GPU architectures \cite{luszczek2017towards}.  We consider both standard PCG and flexible PCG algorithms, and investigate the trade-offs between potentially faster convergence with the additional cost per iteration when using this preconditioning approach. Numerical experiments to illustrate performance are reported for deconovlution and image deblurring problems. To have full control over numerical precision, all half precision computations are simulated in MATLAB using Higham's {\tt chop} function \cite{higham2019simulating}. 

We begin by providing some background and establish notation in Section~\ref{sec:Background}. The basic idea of our preconditioning approach is described in Section~\ref{sec:KronSVD},  and we discuss some issues related to low-precision arithmetic simulations in Section~\ref{sec:LowP}. Numerical results are presented in Section~\ref{sec:result}, and we summarize with concluding remarks in Section~\ref{sec:ConcludingRemarks}.

\section{Background and Notation}
\label{sec:Background}

This section is intended to establish notation used throughout the paper, as well as to provide background material on inverse problems, regularization and preconditioning.

\subsection{Tikhonov Regularization} 
\label{sec:tik}

Because $A$ is severely ill-conditioned, ignoring the noise (even if it is small) in (\ref{invp}) and na\"ively solving $A\bfx=\bfb$ (or $\min\||\bfb - A\bfx\|_2$) typically produces very poor approximations of $\bfx\true$. This can be seen from a singular value decomposition (SVD). 

Specifically, suppose\footnote{To simplify notation, we assume $A$ is an $n \times n$ matrix, but everything we discuss in the paper is applicable to $m \times n$ matrices.}  $A = U\Sigma V^T$ is the SVD of $A \in \R^{n \times n}$, where $U$ and $V$ are orthogonal matrices, and $\Sigma = \,$diag$(\sigma_1, \sigma_2, \ldots, \sigma_n)$ is a diagonal matrix, with the singular values satisfying $\sigma_1 \geq \sigma_2 \geq \cdots \geq \sigma_n \geq 0$. 
For inverse problems, the singular values of $A$ decay rapidly toward zero, so even if we assume $A$ is nonsingular and compute the inverse solution $\bfx\inv = A^{-1}\bfb$, we obtain
$$
\bfx\inv = A^{-1}(A\bfx\true + \bfeps) = \bfx\true + A^{-1}\bfeps  = \bfx\true +\sum_{i=1}^n\frac{\bfu_i^T\bfeps}{\sigma_i}\bfv_i\,,
$$
where $\bfu_i$ is the $i$-th right singular vector of $A$ (i.e., $i$-th column of $U$), and $\bfv_i$ is the $i$-th left singular vector of $A$ (i.e., $i$-th column of $V$).
That is, $\displaystyle \bfx\inv = \bfx\true + $ error, and we observe that noise components $\bfu_i^T\bfeps$ can be greatly magnified when divided by tiny $\sigma_i$ in the error term, $\displaystyle \sum_{i=1}^n\frac{\bfu_i^T\bfeps}{\sigma_i}\bfv_i$. Another important property of inverse problems is that the singular vectors $\bfv_i$ corresponding to large singular values tend to be smooth (low frequency), and $\bfv_i$ corresponding to small singular values tend to highly oscillatory (high frequency). Thus, division by tiny singular values, $\sigma_i$,  magnifies the high frequency components in the error, $\bfv_i$, and these high frequency components dominate any information about $\bfx\true$ in the inverse solution $\bfx\inv$.

There are many regularization approaches that can be used to mitigate noise magnification \cite{CalSom07,EnHaNe00,hansen2010discrete,MuSi12,Vogel02}. In this paper we use Tikhonov regularization, which can be written as 
\begin{equation}\label{eq:Tik}
   \min_\bfx \|\bfb-A\bfx\|^2_2+\lambda^2\|\bfx\|^2_2
    =\min_\bfx \left\|\left[ \begin{array}{c}
\bfb\\
\bfzero\\
\end{array} \right]-\left[ \begin{array}{c}
A\\
\lambda I\\
\end{array} \right]\bfx\right\|^2_2 \, . 
\end{equation}
The regularization parameter $\lambda$ needs to be chosen to balance fit-to-data ($\|\bfb - A \bfx\|_2^2$) with smoothness of the solution ($\|\bfx\|_2^2$); if $\lambda$ is tiny, then we get good fit-to-data with a solution that is highly magnified by noise, and if we choose $\lambda$ too large, then we get a smooth solution that does not provide a very good fit-to-data. It is also interesting to see that if we replace $A$ with its SVD, and compute the least squares solution of (\ref{eq:Tik}), we obtain
\begin{equation}
\label{eq:TikSV}
 \bfx =\sum_{i=1}^n  \phi_i \frac{\bfu_i^T\bfb}{\sigma_i}\bfv_i = \sum_{i=1}^n  \phi_i \frac{\bfu_i^T\bfb\true}{\sigma_i}\bfv_i  + \sum_{i=1}^n  \phi_i \frac{\bfu_i^T\bfeps}{\sigma_i}\bfv_i \, ,
\end{equation}
where $\displaystyle \phi_i = \frac{\sigma_i^2}{\sigma_i^2+\lambda^2}$ acts as a filter; when $\sigma_i$ is large compared to $\lambda$, $\phi_i \approx 1$ allowing low frequency components of the inverse solution (those corresponding to large singular values) to be reconstructed, but when $\sigma_i$ is small compared to $\lambda$, $\phi_i \approx 0$ and the high frequency components in the error term (those corresponding to small singular values) are suppressed. 

\subsection{Finding Regularization Parameters}
\label{sec:RegParams}

Choosing the regularization parameter $\lambda$ is crucial to obtaining a good approximate solution. In this subsection we review a few approaches that can be used, provided we are able to compute the SVD, $A = U\Sigma V^T$.

For experimental purposes, it can be helpful to use test problems where the true solution, $\bfx\true$ is known, and find a parameter that minimizes the relative error; that is, find
\begin{equation}\label{eq:OptTik}
\lambda\opt = \arg\min_\lambda \frac{\|\bfx_\lambda - \bfx\true\|_2}{\|\bfx\true\|_2}\,, \quad \mbox{where} \quad
\bfx_\lambda = \arg\min_\bfx 
\left\|\left[ \begin{array}{c} \bfb\\ \bfzero \end{array} \right] - 
\left[ \begin{array}{c} A \\ \lambda I \end{array} \right] \bfx \right\|^2_2\,.
\end{equation}
Substituting $A = U\Sigma V^T$ into (\ref{eq:OptTik}), it is not difficult to show that $\lambda\opt$ can be found by minimizing a one-dimensional function,
\begin{equation}
\label{eq:OptTikSVD}
\lambda\opt = \arg\min_\lambda \CO(\lambda) = \arg\min_\lambda \sum_{i=1}^n \left(\frac{\sigma_i\hat{b}_i}{\sigma_i^2+\lambda^2}-\hat{x}_i\right)^2\, ,
\end{equation}
where $\hat\bfb = U^T\bfb$ and $\hat\bfx = V^T\bfx\true$. This minimization problem can be found in MATLAB using, for example, the built-in function {\bf fminbnd}, assuming $\sigma_n \leq \lambda \leq \sigma_1$. We remark that this approach is not practical for real problems, but it can be a good reference point (or benchmark) in numerical tests, when comparing with other approaches that do not require knowledge of $\bfx\true$.

The generalized cross validation (GCV) method is a popular approach that does not require $\bfx\true$. The goal of GCV is to predict any missing data from the data set. This means that the difference between an arbitrary $b_i$ left-out of the original data set and the corresponding solution filtered by $\lambda$ should be minimized \cite{chung2008weighted}. This minimization is carried out for all data points in $\bfb$. Although it is not obvious from this description, it can be shown (see, e.g., \cite{hansen2010discrete}) that an estimate of the regularization parameter using the GCV approach can be found by minimizing the one-dimensional function
\begin{equation}
\label{eq:GCVTik}
    \CG(\lambda) = \frac{n\|(I-AA_\lambda^\dag)\bfb\|_2^2}{\left(\text{trace}(I-AA_\lambda^\dag)\right)^2}\, .
\end{equation}
Here $A_\lambda^\dag=(A^TA+\lambda^2I)^{-1}A^T$ and trace means summing up the diagonal elements of a matrix. Using $A = U\Sigma V^T$, it can be shown \cite{hansen2010discrete} that the GCV chosen regularization parameter can be found by the minimization
\begin{equation}
\label{eq:GCVTikSVD}
    \lambda\gcv = \arg\min_\lambda \CG(\lambda) = \arg\min_\lambda \frac{n \displaystyle \sum_{i=1}^n\left(\frac{\hat{b}_i}{\sigma_i^2+\lambda^2}\right)^2}{\left(\displaystyle \sum_{i=1}^n\frac{1}{\sigma_i^2+\lambda^2}\right)^2}\, .
\end{equation}
Note that although equation (\ref{eq:GCVTikSVD}) appears to be complicated, it is a simple one-dimensional minimization problem for $\lambda$, which can be solved, for example, using MATLAB's \textbf{fminbnd} function.

According to \cite{chung2008weighted}, GCV can sometimes choose a $\lambda$ that is too small, and that a weighted version of GCV (wGCV) can overcome this problem:
\begin{equation}
\label{eq:wGCVTik}
     \CG_\omega(\lambda) = \frac{n\|(I-AA_\lambda^\dag)\bfb\|_2^2}{\left(\text{trace}(I-\omega AA_\lambda^\dag)\right)^2}\, .
\end{equation}
Using $A = U\Sigma V^T$, the wGCV chosen regularization parameter (for Tikhonov regularization) can be found by the minimization
\begin{equation}
\label{eq:wGCVTikSVD}
    \lambda\wgcv = \arg\min_\lambda \CG_\omega(\lambda) = 
    \arg\min_\lambda \frac{n \displaystyle \sum_{i=1}^n\left(\frac{\lambda^2\hat{b}_i}{\sigma_i^2+\lambda^2}\right)^2}{\left(\displaystyle \sum_{i=1}^n\frac{(1-\omega)\sigma_i^2+\lambda^2}{\sigma_i^2+\lambda^2}\right)^2}\,.
\end{equation}
When the weight parameter $\omega = 1$ we get the standard GCV function, $\omega > 1$ gives us smoother solutions, and solutions for $\omega < 1$ are less smooth.

Another well-known method for estimating regularization parameters is the discrepancy principle. The main idea is to set an upper bound for the residuals $\|A\bfx_\lambda-\bfb\|_2$, where $\bfx_\lambda$ is the solution computed with regularization parameter $\lambda$. That is, we look for a regularization parameter $\lambda$ such that
\begin{equation}
\label{eq:DSCTik}
    \CD(\lambda) = \|A\bfx_\lambda -\bfb\|_2^2-\epsilon^2 = 0\, ,
\end{equation}
where $\epsilon = \eta\|\bfeps\|$, and $\eta \approx 1$. As with the previous strategies for estimating regularization parameters, if we use $A = U\Sigma V^T$ then the function $\CD(\lambda)$ can be simplified. Specifically, the discrepancy principle reduces to a root finding problem, where we find $\lambda = \lambda\disc$ such that
\begin{equation}
\label{eq:DSCTikSVD}
\CD(\lambda)=\sum_{i=1}^n\left(\frac{\lambda^2\hat{b}_i}{\sigma_i^2+\lambda^2}\right)^2-\epsilon^2=0\, .
\end{equation}
The built-in MATLAB function \textbf{fzero} can be used to find $\lambda\disc$. 

We remark that to use the discrepancy principle, it is necessary to have a good estimate of the noise, $\epsilon = \|\bfeps\|_2$.

\subsection{Iterative approaches for Tikhonov regularization}
\label{sec:PCG}

For large-scale problems, it may not be computationally feasible to compute an SVD unless $A$ has an exploitable structure. For example, if $A$ is a circulant matrix, then fast Fourier transforms (FFTs) can be used to compute a spectral decomposition, and if $A$ is a block circulant matrix with circulant blocks, then 2-dimensional FFTs can be used. Another example occurs if $A$ can be written as a Kronecker product of two smaller matrices; 
we discuss this in more detail in Section~\ref{sec:KronSVD}.

If it is not computationally feasible to compute an SVD of $A$, then iterative methods, such as 
the conjugate gradient (CG) algorithm, are often recommended. 
For least squares problems, CG can be applied to the normal equations, which for the Tikhonov regularization problem (\ref{lsp:Tik})takes the form
\begin{equation*}
    \left[ \begin{array}{cc}
A^T&\lambda I\\
\end{array} \right]\left[ \begin{array}{c}
A\\
\lambda I\\
\end{array} \right]\mathbf{x} = \left[ \begin{array}{cc}
A^T&\lambda I\\
\end{array} \right]\left[ \begin{array}{c}
\mathbf{b}\\
\mathbf{0}\\
\end{array} \right] \, .
\end{equation*}
\begin{equation}
\label{tik:normal2}
(A^TA+\lambda^2I)\mathbf{x}=A^T\mathbf{b}\, .
\end{equation}
Normally, it is recommended to use either CGLS or LSQR implementations \cite{bjorck1990least}. The overall cost of this algorithm depends on the number of iterations and the cost of matrix-vector multiplications involving $A$ and $A^T$, which can be done efficiently if $A$ is sparse or if it has an exploitable structure (e.g., Toeplitz or a sum of Kronecker products).

Preconditioning can be used in the CG algorithm (PCG) to accelerate convergence. The basic idea, for equation (\ref{tik:normal2}), is to construct a matrix $M$ that has three properties. The first is that $M \approx A^TA+\lambda^2 I$, such that the eigenvalues of $M^{-1}(A^TA+\lambda^2I)$ are clustered around a point bounded away from zero (the more eigenvalues that are clustered, and tighter clusters, leads to faster convergence). The second is that solving the linear system $M\mathbf{z}=\mathbf{r}$ is not too expensive, because it is required at each iteration of the PCG algorithm. The third property is that it is not too expensive to construct the approximating matrix $M$. 

In this paper we consider preconditioners constructed from Kronecker product SVD approximations of $A$. We discuss Kroneker products later, but for now it is helpful to understand the basic idea. Suppose were are able to efficiently compute the decomposition $\widehat A = \widehat U \widehat\Sigma \widehat V^T \approx A$, and define 
\begin{align*}
    M &= \widehat A^T \widehat A + \lambda^2I \\
    &= \widehat V \widehat\Sigma^T \widehat U^T \widehat U \widehat\Sigma \widehat V^T + \lambda^2 I \\
    &= \widehat V \widehat\Sigma^T \widehat\Sigma \widehat V^T + \lambda^2 I \\
    &= \widehat V (\widehat\Sigma^T \widehat\Sigma + \lambda^2 I) \widehat V^T\, .
\end{align*}
If we can compute such an approximation, then the solution of the linear system $M\bfz = \bfr$ can be computed as
\begin{equation}
\label{eq:PreconSystem}
\bfz = M^{-1}\bfr = \widehat V\big(D^{-1}\big(\widehat V^T\bfr\big)\big)\,, \quad \mbox{where} \quad D = \widehat\Sigma^T \widehat\Sigma + \lambda^2 I\,.
\end{equation}
That is, solving $M\bfz = \bfr$ requires matrix-vector multiplication with $\widehat V$ and $\widehat V^T$, and multiplication by the diagonal matrix $D^{-1}$. 

\section{Kronecker Product SVD Approximations for Preconditioning}
\label{sec:KronSVD}

In this section we discuss the details of the low precision, Kronecker product SVD approximation preconditioner we propose to use in this work. 

We first begin by supposing $A \in \R^{N \times N}$ where $N = n^2$ can be approximated by a {\em Kronecker product}
$$
  A \approx \widehat{A} = A_r \otimes A_c 
     = \left[ \begin{array}{ccc}
       a^{^{(r)}}_{11} A_c & \cdots & a^{^{(r)}}_{1n} A_c \\[3pt]
       \vdots & & \vdots \\[6pt]
       a^{^{(r)}}_{n1} A_c & \cdots & a^{^{(r)}}_{nn} A_c
     \end{array} \right]
$$
where $A_r$ and $A_c$ are $n \times n$ matrices, and $a^{^{(r)}}_{ij}$ is the $(i,j)$ entry of $A_r$.   There are many convenient properties of Kronecker products that can be exploited for efficient computations \cite{vanloan2000ubiquitous}. Specifically, for general matrices $C, D, E, F \in \R^{n \times n}$ and $\bfy \in \R^{n^2}$, 
\begin{itemize}
\item
$(C \otimes D)(E \otimes F) = CE \otimes DF$.
\item
$(C \otimes D)^T = C^T \otimes D^T$.
\item 
If $C$ and $D$ are nonsingular, then so is $C \otimes D$, and 
$(C \otimes D)^{-1} = C^{-1} \otimes D^{-1}$.
\item
$\bfx = (C \otimes D)\bfy = {\tt vec}(DYC^T)$, where 
$$
\bfy = \mbox{vec}(Y) = \mbox{vec}\left(\left[\begin{array}{cccc} \bfy_1 & \bfy_2 & \cdots & \bfy_n \end{array} \right] \right) =
\left[ \begin{array}{c} \bfy_1 \\ \bfy_2 \\ \vdots \\ \bfy_n \end{array} \right]\,, \quad \bfy_i \in \R^{n}\,.
$$
\end{itemize}
Using these properties for $A \approx A_r \otimes A_c$, with $A_r = U_r \Sigma_r V^T_r$ and $A_c = U_c \Sigma_c V^T_c$, then 
\begin{equation}\label{eq:AKronSVD}
    A \approx \widehat{A} = A_r \otimes A_c = (U_r \Sigma_r V^T_r) \otimes (U_c \Sigma_c V^T_c) = (U_r \otimes U_c)(\Sigma_r \otimes \Sigma_c) (V_r \otimes V_c)^T\, .
\end{equation}
The matrices $U_r \otimes U_c$ and $V_r \otimes V_c$ do not need to 
be formed explicitly, thus required storage for the
decomposition is only $O(N)$.  
(We do not mention anything about $\Sigma_r \otimes \Sigma_c$
because it involves diagonal matrices, and so contributes only
trivially to the computational costs.)
In addition, the decomposition
and computations with $U_r \otimes U_c$ and $V_r \otimes V_c$
require only $O(N^{3/2})$ arithmetic operations.
This is very efficient compared to the standard approaches applied directly
to $A$, which require $O(N^2)$ storage and $O(N^3)$ operations.

A general approach for constructing a Kronecker product approximation was proposed by Van Loan and Pitsianis \cite{vanloan1993approximation}, and an alternative approach for matrices that have a block Toeplitz with Toeplitz block (BTTB) structure, which arise in image deblurring problems, was described in \cite{kamm1998kronecker,kamm2000optimal,nagy2006kronecker,nagy2004kronecker}. An implementation to efficiently construct this approximation for image deblurring problems, which we use for this work, is in the IR Tools MATLAB software package \cite{gazzola2019irtools}. 

To solve linear systems (\ref{eq:PreconSystem}) with the Kronecker product SVD approximation (\ref{eq:AKronSVD}), we exploit the properties of Kronecker products. Specifically, the vector $\bfz$ in equation (\ref{eq:PreconSystem}) can be computed as follows,
\begin{equation}
\label{eq:PreconSolve}
Z = V_c \Big(S \dotmult \big(V_c^T R \,V_r\big) \Big)V_r^T\,, 
\end{equation}
where we use the MATLAB notation $\dotmult$ for element-wise multiplication, $\mbox{vec}(S) = \mbox{diag}(D^{-1})$, $\bfr = \mbox{vec}(R)$, and $\bfz = \mbox{vec}(Z)$.

Another important aspect of our preconditioner is that we store the components of $M$ (i.e., $V_r$, $V_c$ and $S$), and perform the computations in (\ref{eq:PreconSolve}), 
using 16-bit (half) precision arithmetic. Although we do not compute $S$, $V_r$ and $V_c$ in half precision, this should be considered in future work. 
We note that using half precision for preconditioner solves means that, due to roundoff errors, the preconditioner does not remain constant at each iteration, and could lead to a loss of orthogonality of the generated Krylov search directions in PCG \cite{anzt2019adaptive}. A workaround to this problem is to use a flexible PCG algorithm \cite{notay2000flexible}, which essentially replaces how an update parameter of the search direction (standard PCG uses Fletcher-Reeves) with an alternative formula (Polak-Ribie\'re). This workaround requires an additional inner product and additional storage \cite{anzt2019adaptive}.

\section{Low Precision Simulation in MATLAB}
\label{sec:LowP}

Before proceeding to the numerical experiments, we discuss a few important issues related to simulating low precision arithmetic in MATLAB. 
Higham and Pranesh presented a MATLAB function named \textbf{chop} that can be used to simulate lower precision arithmetic by rounding entries of vectors and matrices to the target precision level \cite{higham2019simulating}. It supports various precision formats, for example IEEE half precision, or even user-customized formats. The input of the function is in either double or single precision, and the output is a value at the specified precision; the number is stored in higher precision type with extra bits set to zero. 

When simulating low-precision arithmetic, \textbf{chop} needs to be applied for each floating point operation. For example, the scalar operation
$x+y\times z$
requires two calls to \textbf{chop}:
$$
{\rm chop}(x+{\rm chop}(y\times z)) \, .
$$
When simulating low-precision computations for large-scale problems, there are two important issues to consider:
\begin{itemize}
\item Each call to {\bf chop} is costly.
\item Proper implementations are needed to reduce buildup of roundoff error when replacing built-in computatioanl kernels (e.g., norm) with user defined functions.
\end{itemize}
Consider, for example, the simple computation of an inner product, $s = \bfx^T\bfy$. Since every floating point operation needs to be properly rounded, one might consider the implementation:
\begin{itemize}
\item[]
{\tt s = 0;} \\
{\tt for i = 1:n} \\
\hspace*{12pt} {\tt s = chop(s + chop(x(i)*y(i)))} \\
{\tt end}
\end{itemize}
This requires $2n$ calls to {\bf chop}. We can reduce this to $n+1$ calls by exploiting vectorization:
\begin{itemize}
\item[]
{\tt z = chop(x\ttdotmult y);} \\
{\tt s = 0;} \\
{\tt for i = 1:n} \\
\hspace*{12pt} {\tt s = chop(s + z)} \\
{\tt end}
\end{itemize}
However, it is well-known that the above approach to compute the sum of floating point numbers can lead to a buildup of rounding errors.  Alternative summation approaches include various ordering schemes, blocking, Kahan's compensated summation, and pairwise summation \cite{higham1993accuracy}. In our implementations, we use pairwise summation because, in addition to having good rounding error properties, it also significantly reduces the number of calls to {\bf chop}.

To visualize the pairwise summation approach to compute $s = \bfx^T\bfy$, first compute the element-wise product $\bfz = {\tt chop}(\bfx \dotmult \bfy)$, so that only one vectorized call to \textbf{chop} is needed to compute all entries in $\bfz$. Pairwise summation then proceeds with a divide-and-conquer approach to sum all of the entries in $\bfz$, two at a time. This can be illustrated with a small example, $\bfz \in \R^8$, noticing that we can do pairwise summation as a vectorized operation, followed by a single vectorized call to {\bf chop} at each stage:

\begin{eqnarray*}
\bfz^T = \bfz^{(0)T} & = & 
\left[ \begin{array}{cccccccc} 
\bfz^{(0)}_1 & \bfz^{(0)}_2 & \bfz^{(0)}_3 & \bfz^{(0)}_4 & \bfz^{(0)}_5 & \bfz^{(0)}_6 & \bfz^{(0)}_7 & \bfz^{(0)}_8 
\end{array} \right] \\[-12pt]
& &
\begin{array}{cccc}
\; \underbrace{\hspace*{42pt}}_{\text{add}} & \; \underbrace{\hspace*{42pt}}_{\text{add}} & \; \underbrace{\hspace*{42pt}}_{\text{add}} & \; \underbrace{\hspace*{42pt}}_{\text{add}} \end{array} \\
& & \hspace*{20pt} \underbrace{\hspace*{182pt}}_{\text{vectorized chop}} \\
& = & 
\left[ \begin{array}{cccc} 
\hspace*{14pt} \bfz^{(1)}_1 \hspace*{14pt} & \hspace*{14pt} \bfz^{(1)}_2 \hspace*{14pt} & \hspace*{14pt} \bfz^{(1)}_3 \hspace*{14pt} & \hspace*{14pt} \bfz^{(1)}_4 \hspace*{14pt}
\end{array} \right] \\[-12pt]
& &
\begin{array}{cc}
\hspace{16pt} \underbrace{\hspace*{78pt}}_{\text{add}} & \hspace*{24pt} \underbrace{\hspace*{78pt}}_{\text{add}} \end{array} \\
& & \hspace*{50pt} \underbrace{\hspace*{130pt}}_{\text{vectorized chop}} \\
& = & 
\left[ \begin{array}{cc} 
\hspace*{46pt} \bfz^{(2)}_1 \hspace*{41pt} & \hspace*{41pt} \bfz^{(2)}_2 \hspace*{40pt}
\end{array} \right] \\[-6pt]
& &
\hspace*{54pt} \underbrace{\hspace*{130pt}}_{\text{add \& final chop}} 
\end{eqnarray*} 
With this approach, a MATLAB implementation of a low-precision inner product $\bfx^T\bfy$, with $\bfx, \bfy \in \R^n$ and $n = 2^p$, could look like:
\begin{itemize}
\item[]
{\tt z = chop(x\ttdotmult  y)';} \\
{\tt for k = 1:p} \\
\hspace*{12pt} {\tt K = kron(speye(n/2),[1;1]);} \\
\hspace*{12pt} {\tt z = chop(z*K);} \\
\hspace*{12pt} {\tt n = n/2;} \\
{\tt end}
\end{itemize}
Some comments about the above implementation:
\begin{itemize}
\item
The statement {\tt K = kron(speye(n/2),[1;1]);} constructs a sparse representation of the matrix
$$
K = I \otimes \left[\begin{array}{c} 1 \\ 1 \end{array} \right] = 
\left[\begin{array}{ccc} 
1 & & \\ 1 & & \\ & 1 & \\ & 1 & \\ & & \ddots
\end{array}\right]
$$
\item
The statement {\tt z*K} computes pair-wise sums of entries in (the row vector) {\tt z}, and resizes to half the length. For example, using notation in the previous diagram, 
$\bfz^{(1)T} = \bfz^{(0)T}K$.
\item
When completed, {\tt z} has one entry, the inner product $\bfx^T\bfy$.
\item
This divide-and-conquer approach requires $\log_2(n)+1$ calls to {\tt chop}, which is significantly less than the $2n$ or even $n+1$ calls required by the previously discribed direct summation approaches.
\item
The code can easily be modified for the case when $n$ is not a power of 2.
\item
An advantage of this approach is that it uses a natural blocking that reduces the buildup of roundoff errors when $n$ is large \cite{higham1993accuracy}. 
\end{itemize}

This pairwise summation procedure can be used for matrix-vector multiplication. Specifically, if computing $A\bfv$, we can first scale the columns of $A$ with entries of $\bfv$ and call {\bf chop} once using vectorization (e.g., {\tt C = chop(A\ttdotmult v')}), and then use pairwise summation to add the columns of {\tt C}, which again requires only $\log(n)$ additional calls to {\bf chop}. Similarly, matrix-matrix multiplication can be performed using pairwise summation of outer products, and is the approach we use for our low precision preconditioner solve computation (\ref{eq:PreconSolve}).

\section{Numerical Results}
\label{sec:result}

In this section we provide some numerical experiments to illustrate the effectiveness of using a half-precision Kronecker product SVD approximation as a preconditioner for Tikhonov regularized least squares problems. We use Higham's {\bf chop} package \cite{higham2019simulating} to simulate half-precision for the preconditioner solves, and the IR Tools package \cite{gazzola2019irtools} to generate image deblurring test problems. Specifically, we consider {\em Gaussian} (often used for image deblurring test problems), {\em speckle} (used to simulate blurred observations in telescope imaging of space objects), {\em defocus} (i.e., out-of-focus lens imaging), and {\em shake} (simulating random camera movements) blurs. 

The IR Tools package allows to easily generate $A$ (stored in a structured format that allows for fast matrix-vector multiplications), $\bfx\true$ and $\bfb\true = A\bfx\true$. For each test problem we add noise $\bfeps$ to $\bfb\true$ to obtain the simulated measured data,
$$
\bfb = A \bfx\true + \bfeps\,,
$$
and for the various tests we report the amount of noise as
$$
\mbox{noise level} = \frac{\|\bfeps\|_2}{\|\bfb\true\|_2}\,.
$$
For each test problem, we use $128 \times 128$ images, which means $\bfx\true, \bfb \in \R^{16384}$ and $A \in \R^{16384 \times 16384}$.

For each blur, it is possible to efficiently decompose $A$ into a sum of Kronecker products \cite{kamm1998kronecker,kamm2000optimal},
\begin{equation}
\label{eq:KronSumDecomp}
A = \sum_{i=1}^k A_r^{(i)} \otimes A_c^{(i)}\,
\end{equation}
where 
$$
A_r^{(1)} \otimes A_c^{(1)} = \arg\min_{A_r,A_c} \|A - A_r \otimes A_c\|_F\,.
$$
We can therefore use the Frobenius relative error, $\|A - \widehat A\|_F/\|A\|_F$ as a rough measure of the quality of the preconditioner's approximation of $A$, where $\widehat A$ is given in equation (\ref{eq:AKronSVD}). These norms, for each of the test problems, is given in Table~\ref{tab:PreconApprox}. For the Gaussian blur, $A$ is exactly decomposed into a single Kronecker product, so the relative error between $A$ and $\widehat A$ is due only to the half-precision approximation. 

\begin{table}[!htb]
\caption{Preconditioner approximations for various test problems.}
\label{tab:PreconApprox}
\begin{center}
\begin{tabular}{c|c|c} \hline\hline
 & \# of terms in Kron sum & \\
 Blur type & (i.e., $k$ in equation (\ref{eq:KronSumDecomp})) & $\|A - \widehat A\|_F/\|A\|_F$ \\ \hline\hline
 Gauss & \;\;1 & $7.9195\times10^{-8}$ \\
 Defocus & \;\;6 & $3.1006\times 10^{-1}$ \\
 Speckle & 18 & $4.3102\times 10^{-1}$ \\
 Shake & 10 & $6.2834\times 10^{-1}$ \\
 Motion & \;\;9 & $7.9212\times 10^{-1}$ \\ \hline\hline
 \end{tabular}
 \end{center}
 \end{table}

\subsection{Comparison between Different Parameter Choice Methods} \label{parachoice}

\begin{figure}[!htb]
\begin{center}
\begin{tabular}{cc}
 \includegraphics[width= 0.49 \linewidth]{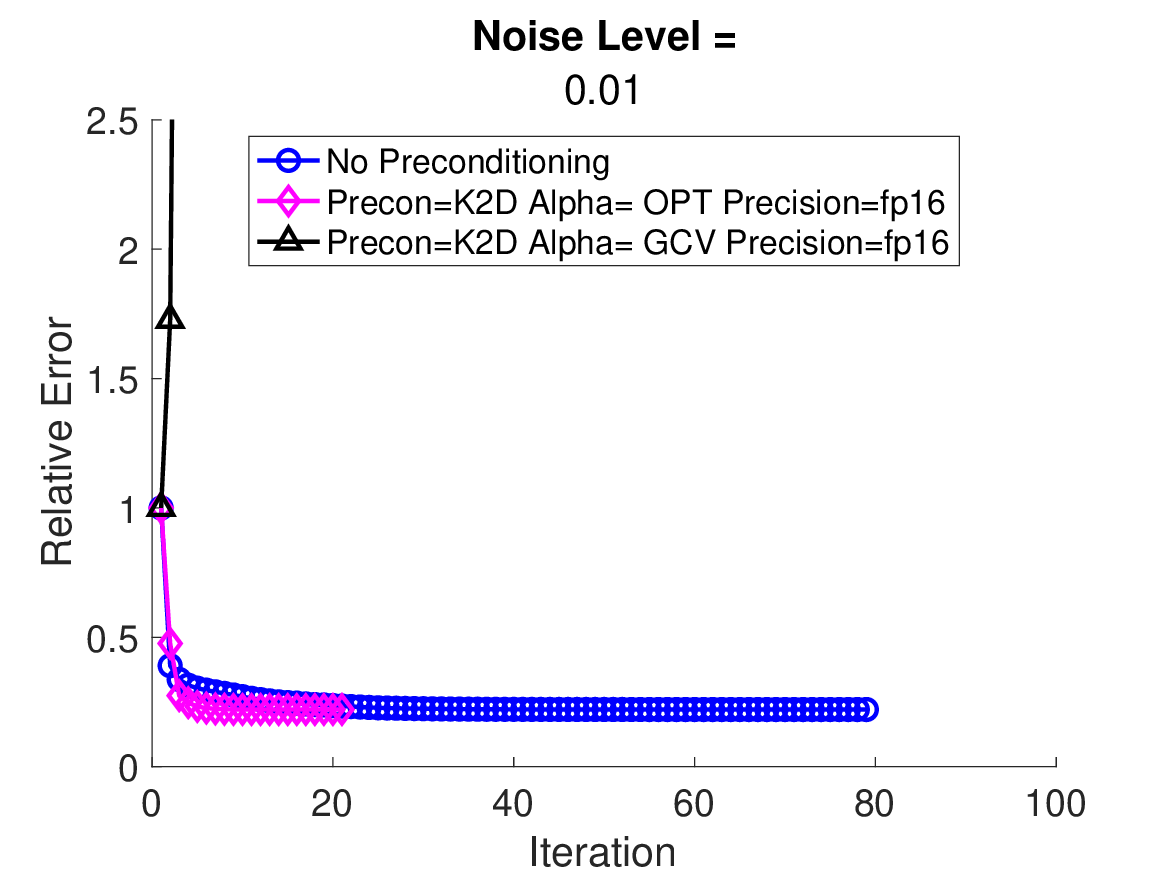} &
 \includegraphics[width= 0.49 \linewidth]{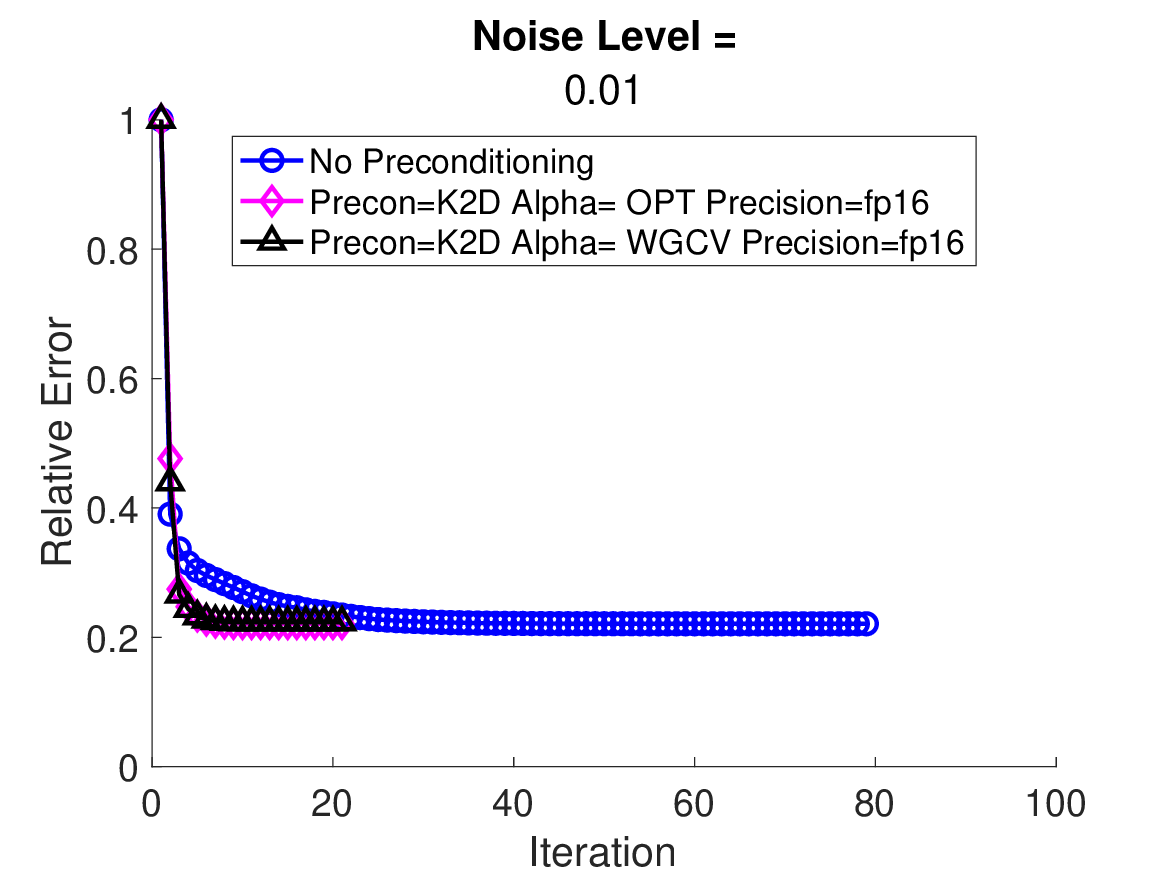} \\
 {(a)} & (b)
 \end{tabular}
 \end{center}
 \caption{A defocus problem of size 128 and 1\% noise (a) with GCV; (b) with weight = $3$ for wGCV}
 \label{gcv_128_defocus_0.01}
\end{figure}

We first consider the issue of determining regularization parameters. The techniques discussed in Section \ref{sec:RegParams} require the SVD of $A$, which is very expensive to compute for large-scale problems. However, because we are using a half-precision Kronecker product SVD approximation for the preconditioner, we will use that to compute estimates of the regularization parameter. While the parameter given by GCV works well for the Gaussian blur type, $\lambda$ tended to be too small for other blur types, such as defocus and speckle, and the relative errors blew up at the second iteration. Figure \ref{gcv_128_defocus_0.01}(a) was a typical example of the instant degeneration.

\begin{figure}[!htb]
\begin{center}
\begin{tabular}{ccc}
 \includegraphics[width= 0.32 \linewidth]{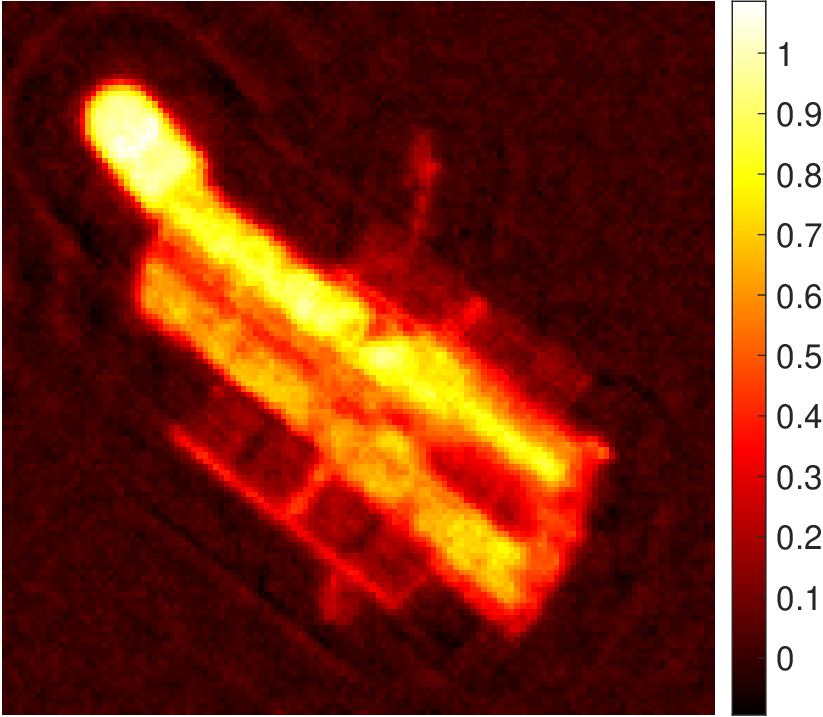} &
 \includegraphics[width= 0.32 \linewidth]{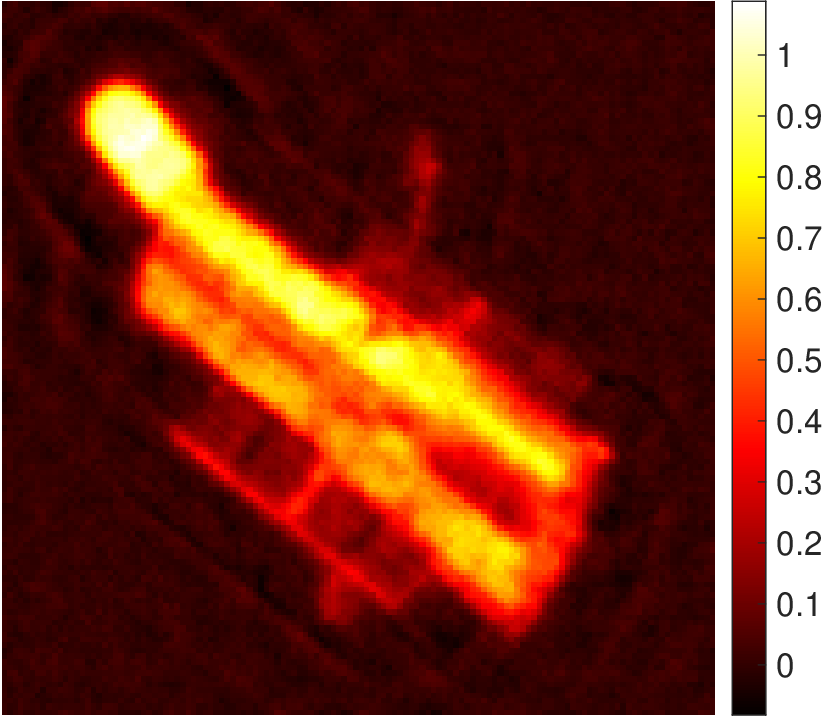}
 &
 \includegraphics[width= 0.32 \linewidth]{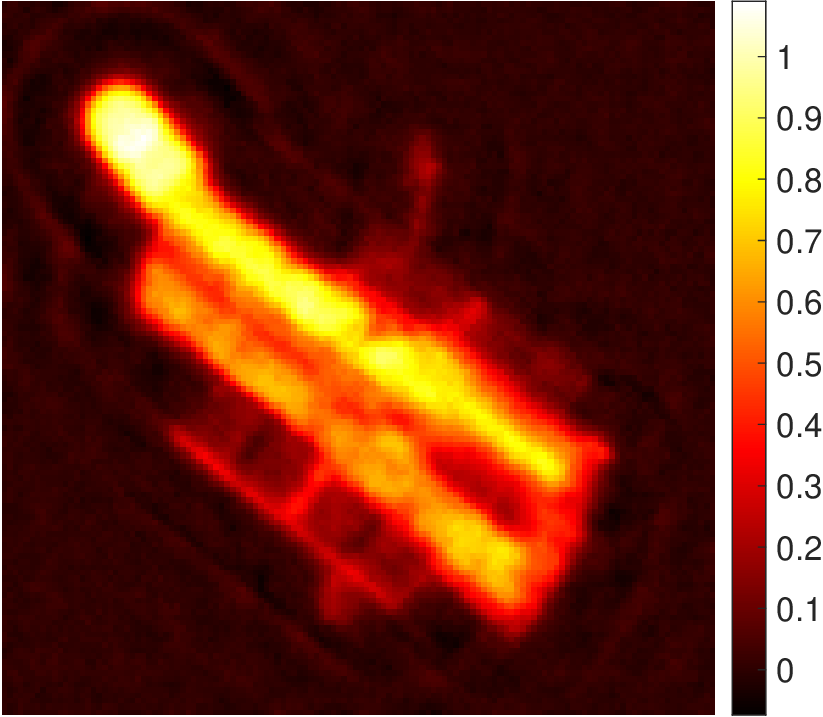}
 \\
 {(a)} & (b) & (c)
 \end{tabular}
 \end{center}
 \caption{Result for a defocus problem of size 128 and 1\% noise with wGCV (a) with weight = $3$; (b) with weight = $5$; (c) with weight = $8$}
 \label{wgcv_128_defocus_0.01}
\end{figure}

Since GCV chose $\lambda$ too small and under-smoothed the solutions, we instead used wGCV and set the weight $\omega >1$. For the same problem as in Figure \ref{gcv_128_defocus_0.01}(a), after applying wGCV with $\omega = 3$, we obtained a solution that was comparable to the one using $\lambda_{opt}$. The vertical axis was rescaled in Figure \ref{gcv_128_defocus_0.01}(b) so that the difference in convergence rate is presented more clearly. The algorithm was not sensitive to the weight: good results could still be obtained if we set $\omega = 5$ or $8$, and the differences in relative errors were negligible, as demonstrated in Figure \ref{wgcv_128_defocus_0.01}.

\begin{figure}[!htb]
\begin{center}
\begin{tabular}{cc}
 \includegraphics[width= 0.49 \linewidth]{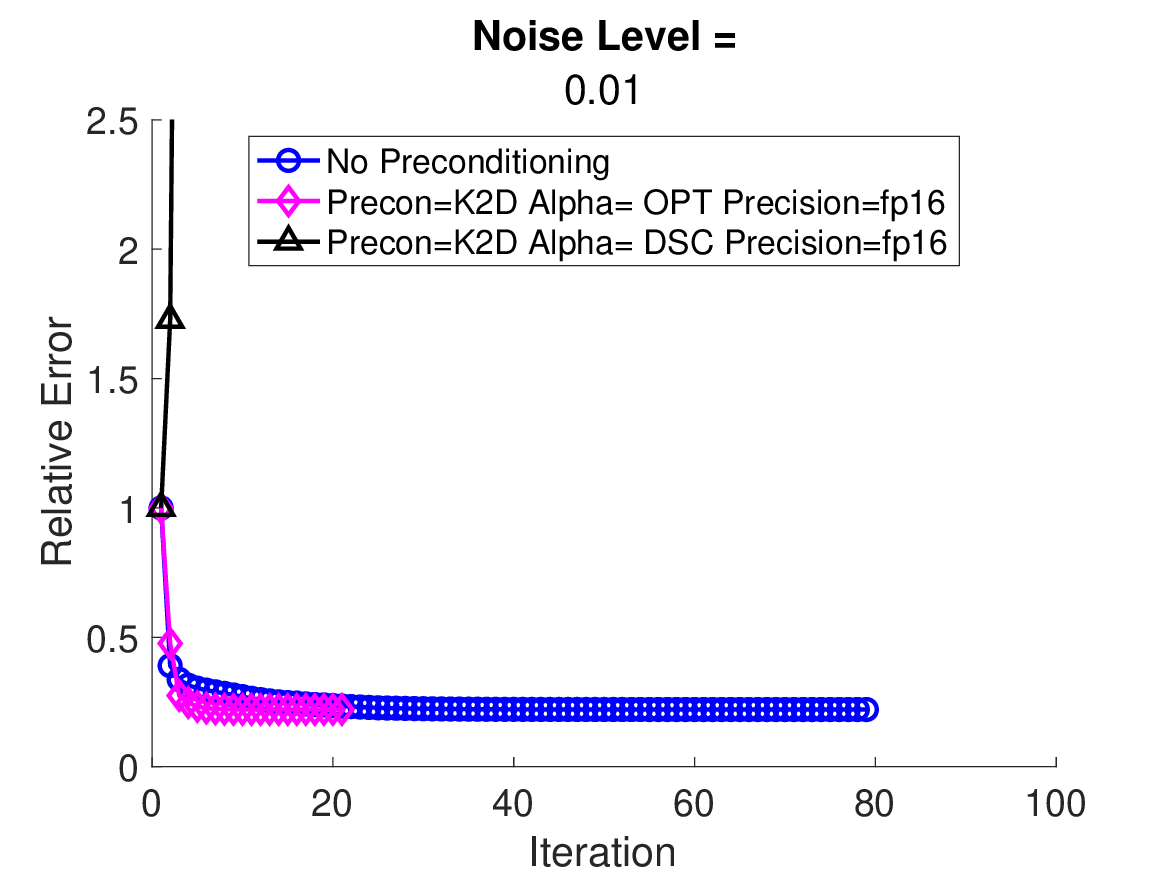} &
 \includegraphics[width= 0.49 \linewidth]{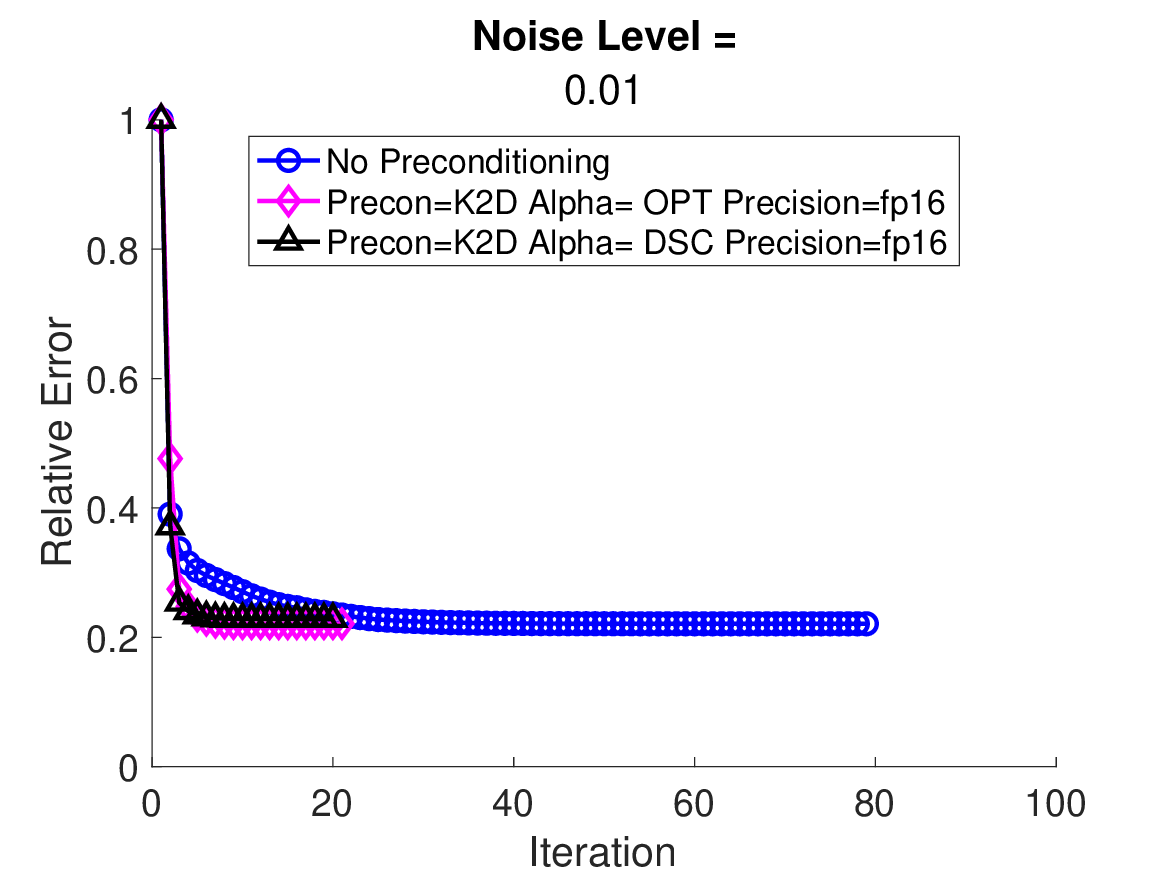} \\
 {(a)} & (b)
 \end{tabular}
 \end{center}
 \caption{A defocus problem of size 128 and 1\% noise (a) with the discrepancy principle; (b) with $\omega = 2$ for the modified discrepancy principle}
 \label{dsc_128_defocus_0.01}
\end{figure}

\begin{figure}[!htb]
\begin{center}
\begin{tabular}{ccc}
 \includegraphics[width= 0.32 \linewidth]{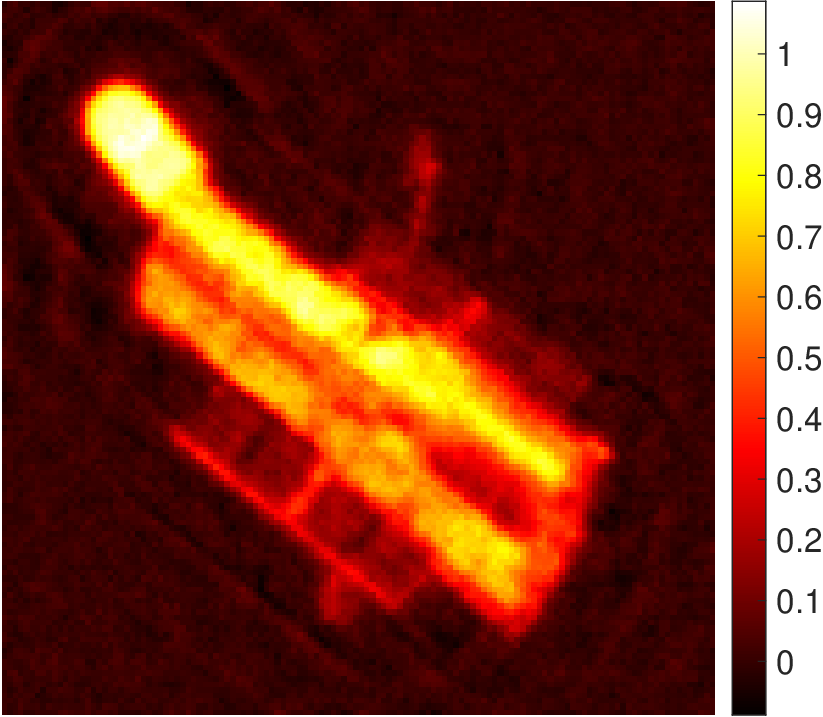} 
 &
 \includegraphics[width= 0.32 \linewidth]{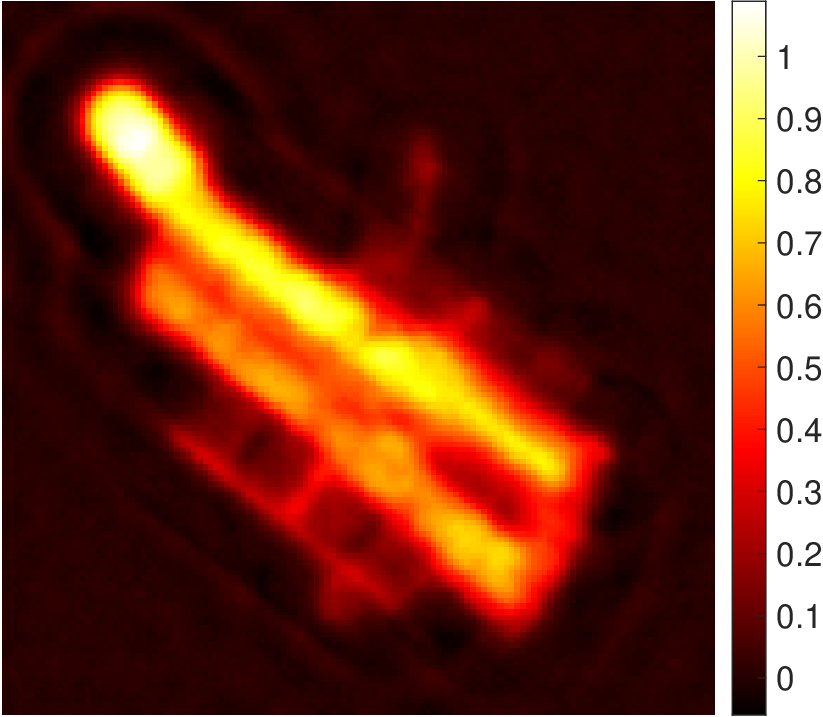}
 &
 \includegraphics[width= 0.32 \linewidth]{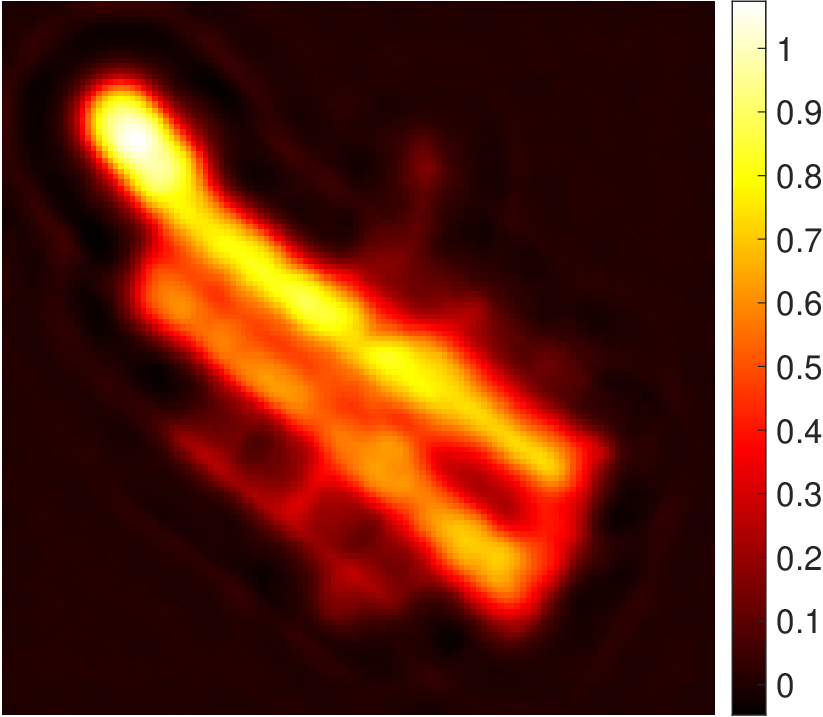}
 \\
 {(a)} & (b) & (c)
 \end{tabular}
 \end{center}
 \caption{Result for a defocus problem of size 128 and 1\% noise with the discrepancy principle (a) with $\eta = 2$; (b) with $\eta = 3$; (c) with $\eta = 5$}
 \label{dsc_128_defocus_0.01_result}
\end{figure}

As for the $\|\bfeps\|_2$-based discrepancy principle, this did not work well with our approximate SVD preconditioner unless we used a large value for $\eta$. 
Figure \ref{dsc_128_defocus_0.01}(b) demonstrated that the discrepancy principle with $\eta = 2$ resulted in a $\lambda$ that was comparable to the performance of $\lambda_{opt}$ for the same problem in Figures \ref{gcv_128_defocus_0.01} and \ref{dsc_128_defocus_0.01}(a). Here the vertical axis was also rescaled for better presentation. However, the discrepancy principle was more sensitive to $\eta$ than the weight parameter in wGCV; there were obvious visible differences in the relative errors and the resulting reconstructed images if $\eta$ was not chosen optimally. Figure \ref{dsc_128_defocus_0.01_result} showed the results for the same defocus problem with different values of $\eta$, and the quality of the solution deteriorated quickly as $\omega$ deviates from its optimal value.

 Moreover, the optimal $\eta$ varied a lot across different problems for the discrepancy principle. Whereas a suitable weight $\omega$ usually ranged from $2$ to $8$ for wGCV, the optimal $\eta$ could be as large as $15$ for certain problem sizes, blur types, and noise levels in the case of the discrepancy principle. A more systematic study of choosing appropriate weight parameters is desired, both for wGCV and the discrepancy principle.

\subsection{Accuracy and Work Unit Evaluation}\label{workuniteval}

\begin{figure}[!htb]
\begin{center}
\begin{tabular}{cc}
 \includegraphics[width= 0.49 \linewidth]{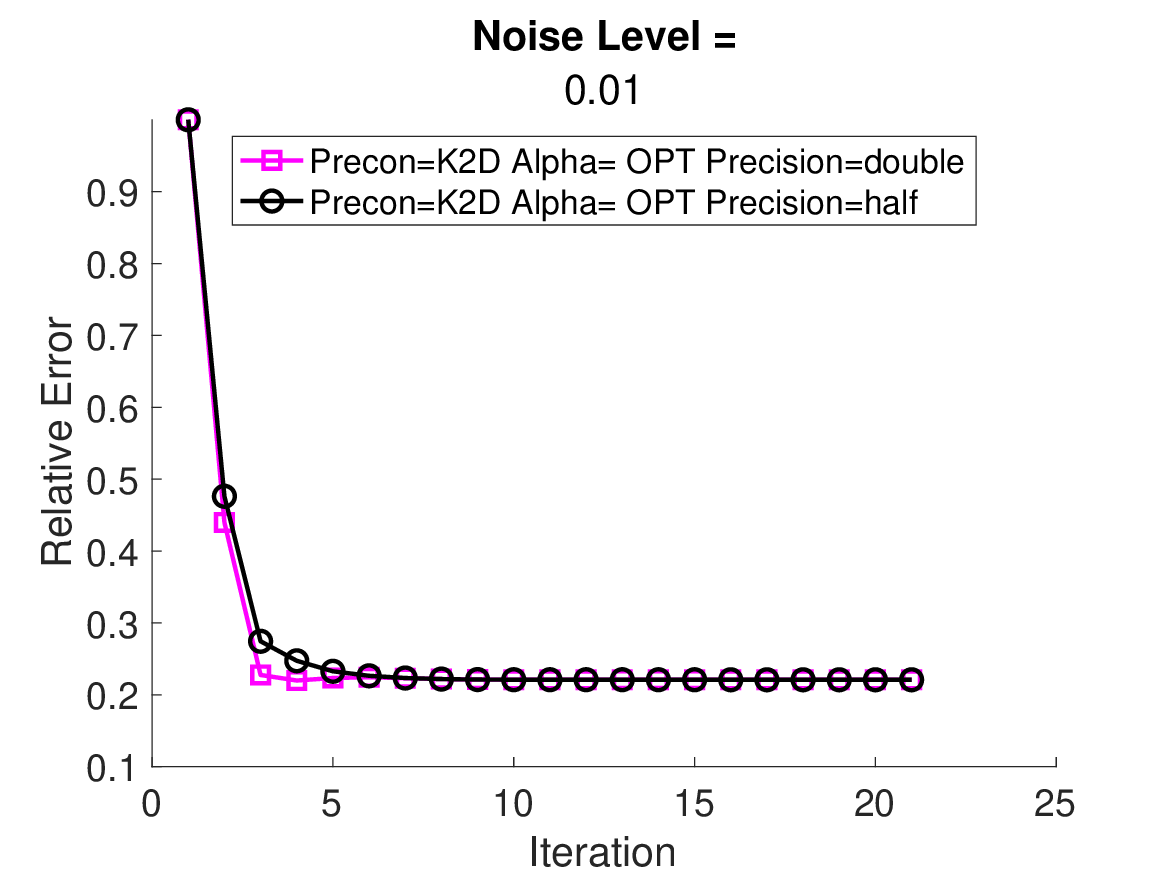} &
 \includegraphics[width= 0.49 \linewidth]{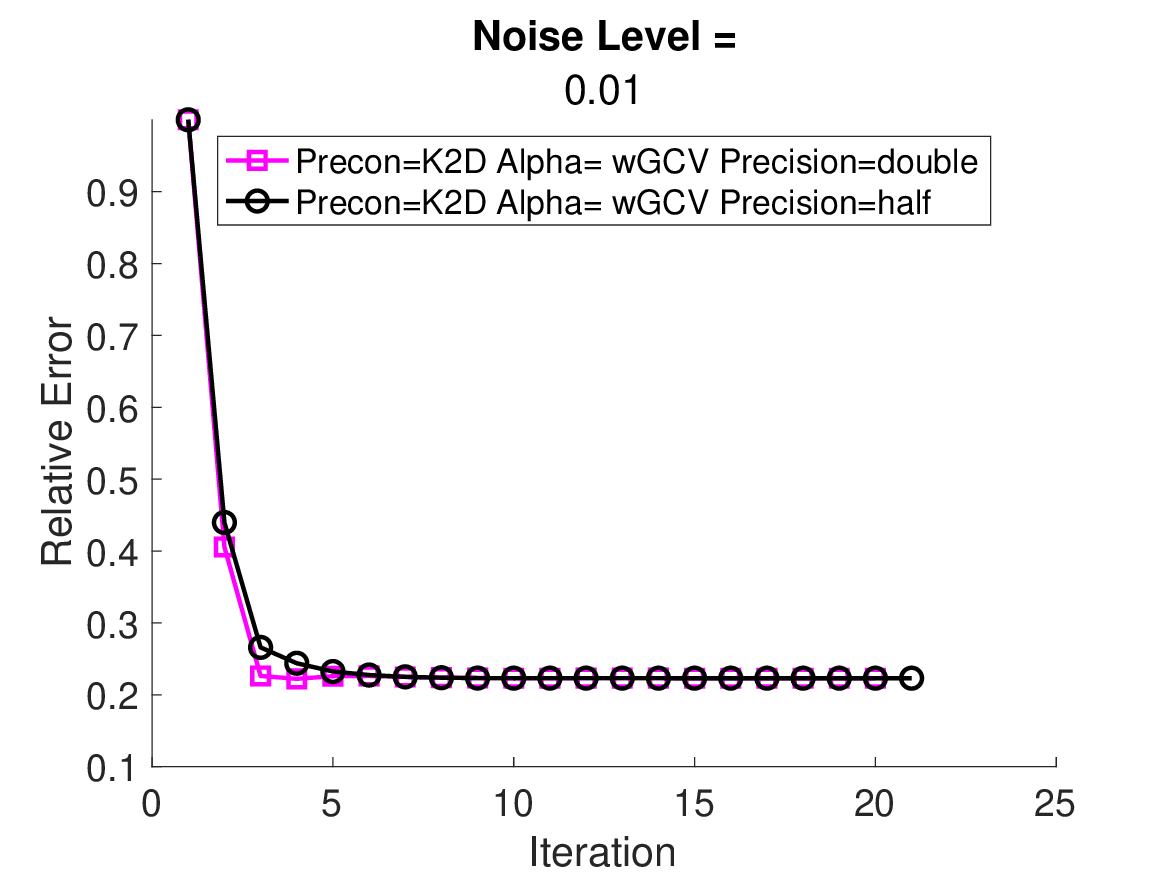} \\
 {(a)} & (b)
 \end{tabular}
 \end{center}
 \caption{Comparison between double and half precision preconditioning for a defocus problem of size 128 and 1\% noise (a) with the optimal method; (b) with wGCV}
 \label{defocus_128_0.01_compare}
\end{figure}

It is certain that reducing the preconditioning to low precision can achieve faster computational speed, but the potential loss in accuracy is also to be considered. Based on the numerical experiments, the difference in relative errors is negligible. For the same test problem as in Section \ref{parachoice}, we compared the relative errors from algorithms that used double-precision preconditioning and half-precision preconditioning, respectively, as shown in Figure \ref{defocus_128_0.01_compare}. While the double precision preconditioning achieved smaller relative errors before convergence, both precision levels converged to the same relative errors in approimately the same number of iterations; it is difficult to tell the two lines apart in the plot. Therefore, although a small increase in the relative errors is incurred when switching to lower precision in the convergence process, the precision level of the preconditioning did not affect the relative errors that the algorithm could reach.

As previously mentioned, using half precision for preconditioner solves means that, due to roundoff errors, the preconditioner does not remain constant at each iteration, and could lead to a loss of orthogonality of the generated Krylov search directions in PCG \cite{anzt2019adaptive}. We tested a flexible PCG method \cite{notay2000flexible}, which requires an additional inner product and additional storage \cite{anzt2019adaptive}, but did not see any difference in convergence behavior. This is probably due to the fact that we do not need many iterations, so loss of orthogonality is not significant.

To get a rough comparison of the additional cost of preconditioning, and to determine if the extra cost is worth preconditioning, we note that 
for each PCG iteration, the most expensive computations are matrix-vector multiplications with $A$ and $A^T$, and the preconditioner solve with $M$. 
Since we are simulating half-precision in software, a direct timing comparison between the preconditioned and unpreconditioned algorithms would not be meaningful.
However, we can get a rough estimate if we assume that each preconditioner solve (in half-precision) is at most $\frac{1}{4}$ the cost of a double-precision matrix-vector multiplication with $A$ or $A^T$. This will be true if we exploit Kronecker product structure, or if we can use fast Fourier transforms (FFTs), which is possible for certain spatially invariant image deblurring problems. 

Let $t$ be the time required for a double-precision matrix-vector multiplication, and let $m_N$ and $m_P$ denote the total number of iterations needed, respectively, for the unpreconditioned and preconditioned algorithms. Then in order for the extra cost of preconditioning to result in a more efficient algorithm, we need
$$
\left(2t + \frac{1}{4}t\right)m_P < 2tm_N \quad \Rightarrow \quad m_P < \frac{8}{9}m_N\,.
$$
In the next section, we see that in most cases we achieve this requirement.

\subsection{Additional Experiments}
So far, we have only tested the algorithm on the defocus problem with 1\% noise. In order to fully understand the advantages and pitfalls of our preconditioning approach, we further applied it to different blur types and noise levels. In Section~\ref{parachoice}, we investigated several parameter choice methods and concluded that both wGCV and the discrepancy principle could be effective with suitably chosen $\omega$ (for wGCV) and $\eta$ (for the discrepancy principle). Therefore, we directly used $\lambda_{opt}$ in this section to assess the maximum potential performance of our algorithm. However, we emphasize that this $\lambda_{opt}$ is computed using the Kronecker product approximate SVD, and not the SVD of the matrix $A$.

\begin{figure}[!htb]
\begin{center}
\begin{tabular}{cc}
 \includegraphics[width= 0.49 \linewidth]{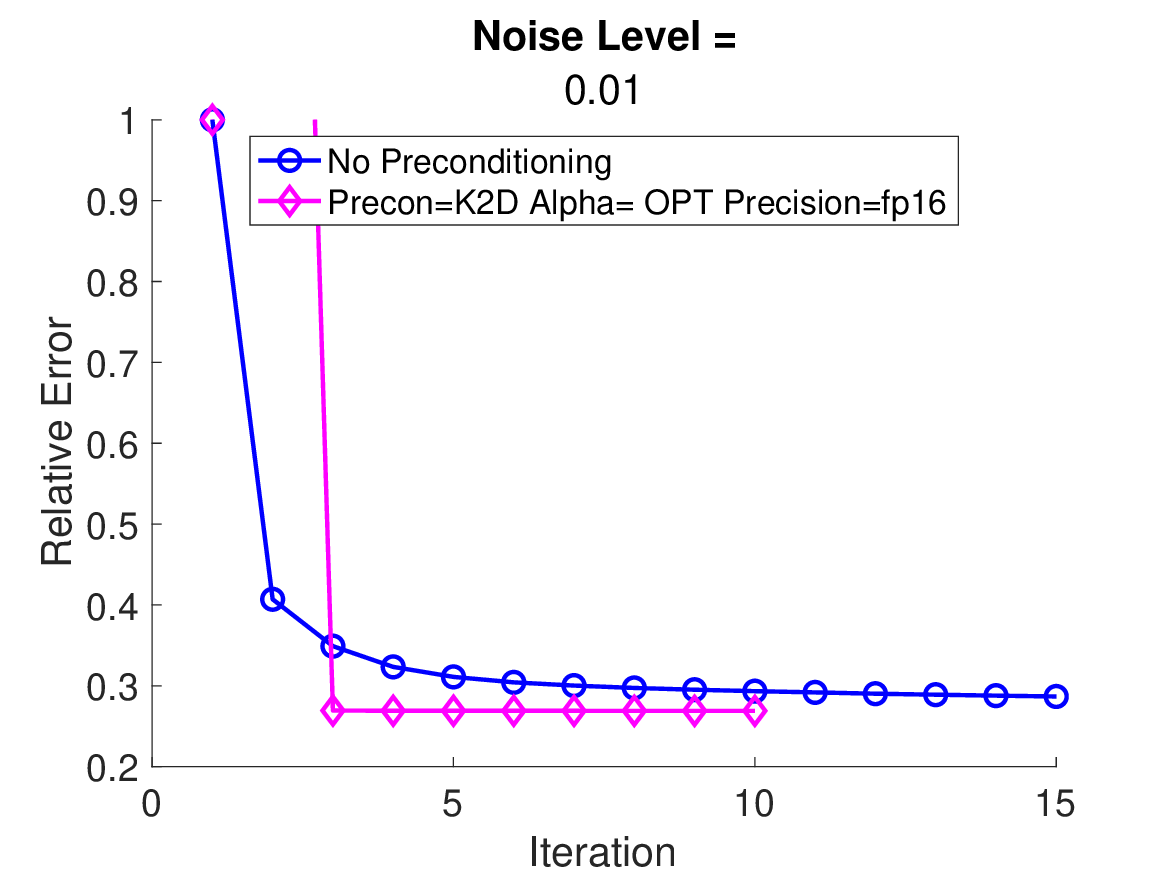} &
 \includegraphics[width= 0.49 \linewidth]{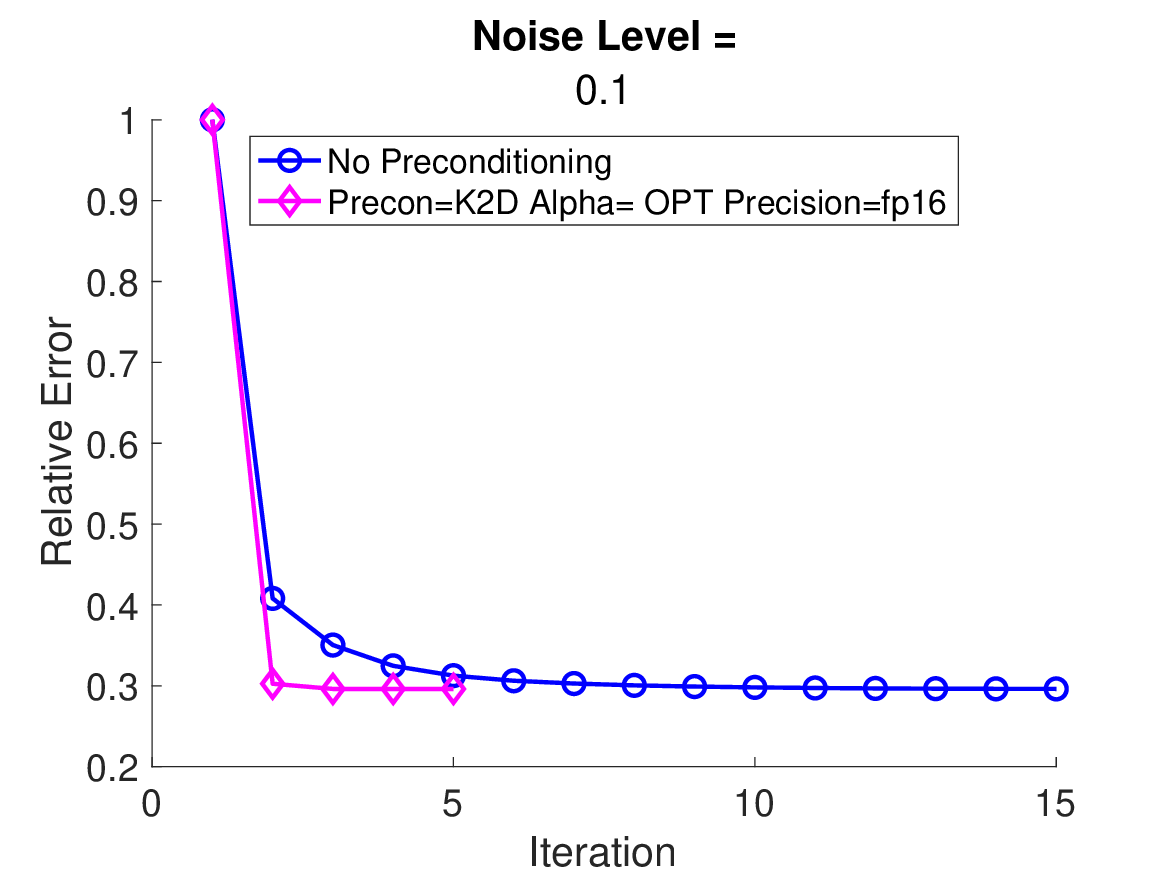} \\
 {(a)} & (b)
 \end{tabular}
 \end{center}
 \caption{A Gaussian problem of size 128 (a) with 1\% noise; (b) with 10\% noise}
 \label{gauss_128}
\end{figure}

First consider the Gaussian blur, where $A$ can be written as a Kronecker product, so the only approximation of the preconditioner comes from using half-precision. As demonstrated by Figure \ref{gauss_128}, in this ideal situation, our algorithm converges at the second iteration, a significant improvement on the no-preconditioning case. In Figure \ref{gauss_128}(a), there is a jump of the relative error in the second iteration. Because we are using the same regularization parameter for both the preconditioner solve and the inverse problem, the parameter is too small for the preconditioner solve in the first few iterations. We could use a larger parameter for the preconditioner solve and eliminate the jump, but since the algorithm quickly recovers from it, we decided to use the same parameter for both scenarios. 

\begin{figure}[!htb]
\begin{center}
\begin{tabular}{cc}
 \includegraphics[width= 0.49 \linewidth]{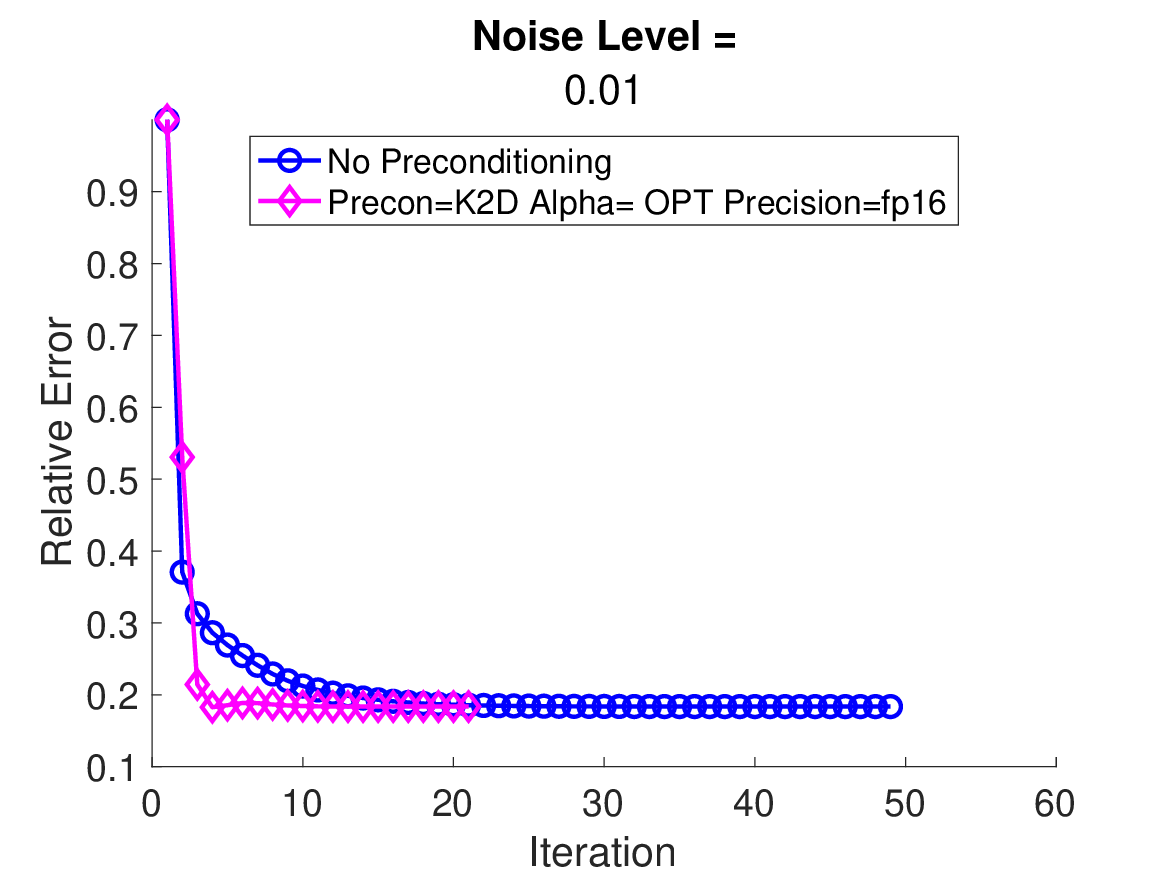} &
 \includegraphics[width= 0.49 \linewidth]{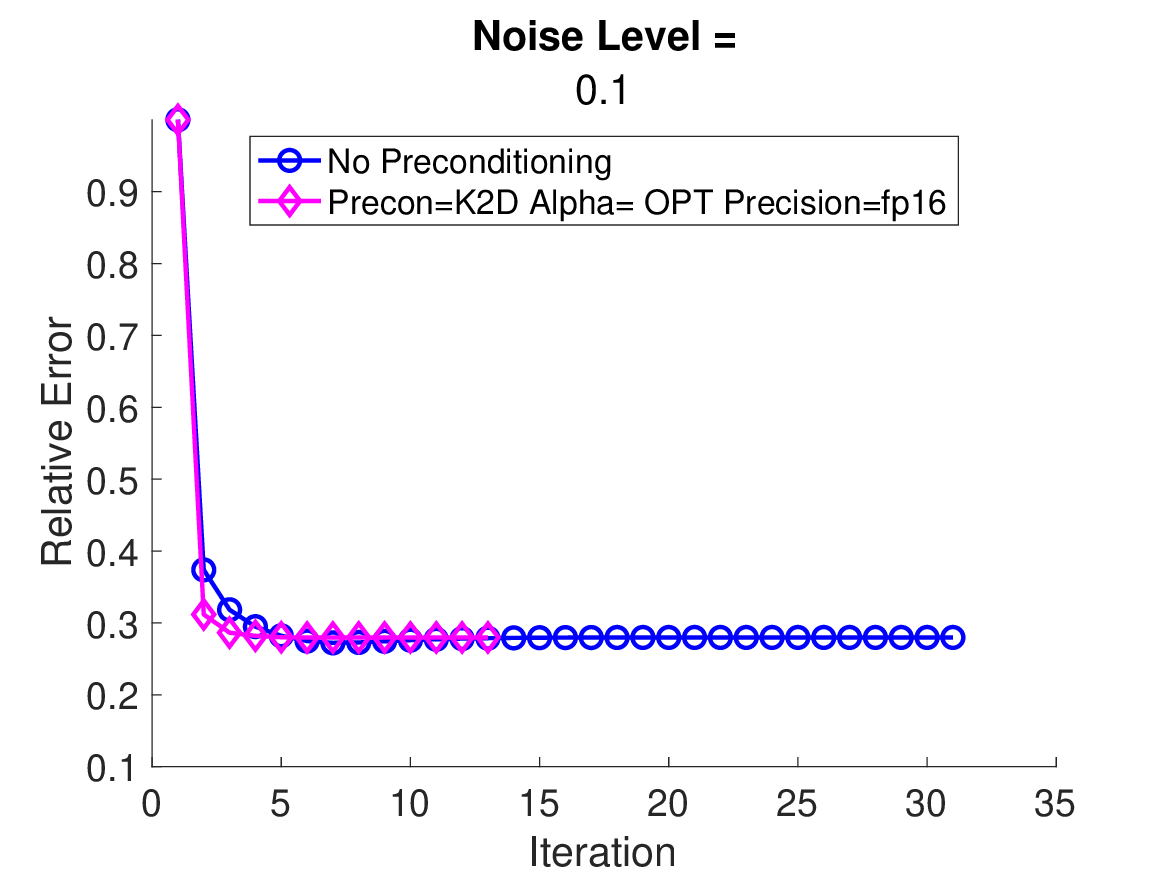} \\
 {(a)} & (b)
 \end{tabular}
 \end{center}
 \caption{A shake problem of size 128 (a) with 1\% noise; (b) with 10\% noise}
 \label{shake_128}
\end{figure}

As its name implies, the shake blur applies random, shaking motion to each pixel \cite{gazzola2019irtools}. As summarized in Table~\ref{tab:PreconApprox}, in this case $A$ is a sum of $10$ Kronecker product terms. For this type of blur, our preconditioning algorithm converges much faster than the no-preconditioning scenario when the noise level is low ($1\%$ in the test case of Figure \ref{shake_128}(a)), but for a high noise level ($10\%$ in the test case of Figure \ref{shake_128}(b)), there is no obvious advantage in preconditioning. One reason behind the lack of improvement from the preconditioning is that we can get little information about the true solution because of the high noise, and therefore the preconditioner is not so useful.

The motion blur creates a linear motion with a $45$-degree angle \cite{gazzola2019irtools}. As summarized in Table~\ref{tab:PreconApprox}, in this case $A$ is a sum of $9$ Kronecker product terms. However, the performance of the proconditioner on this blur type is very different from the shake blur: for both noise levels, the preconditioning provides us with no advantage. We are uncertain of the reason behind such disparity, it is possible that different regularization parameters could have an affect on the convergence behavior.

\begin{figure}[!htb]
\begin{center}
\begin{tabular}{cc}
 \includegraphics[width= 0.49 \linewidth]{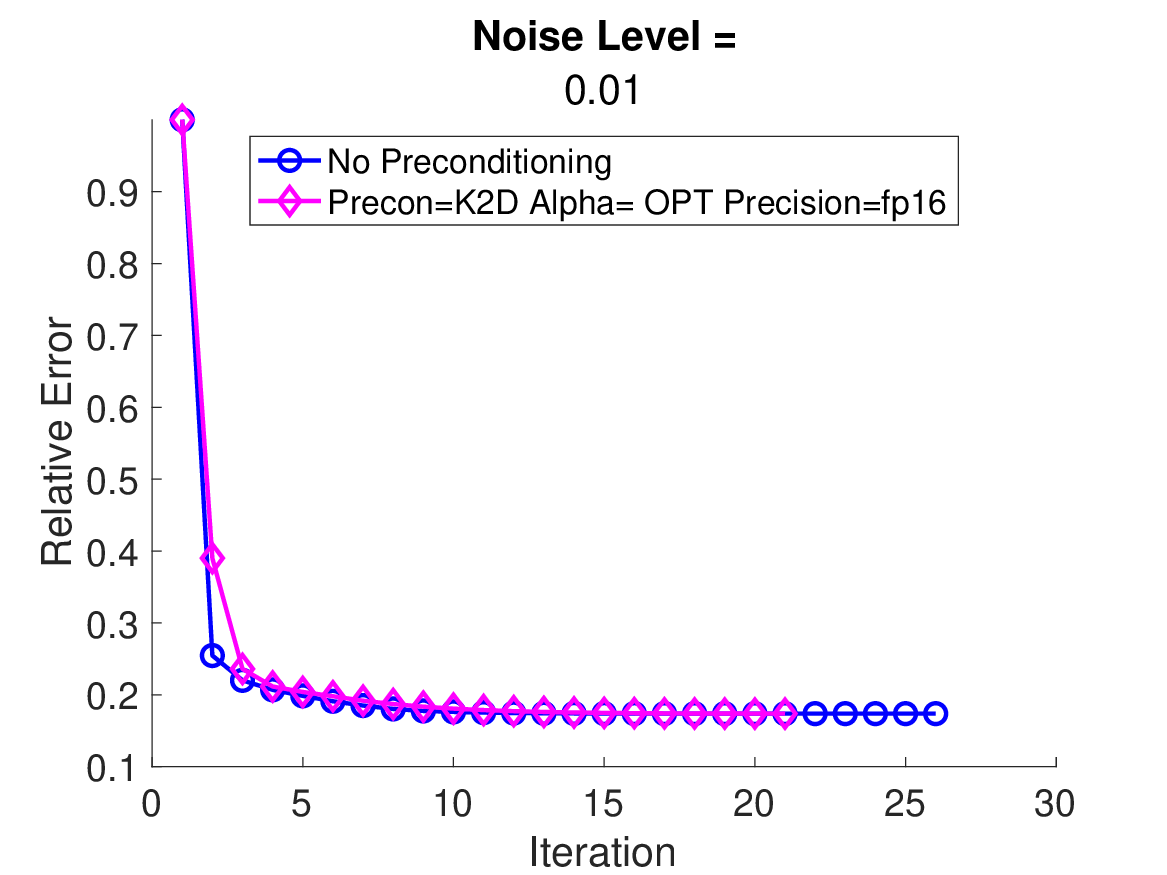} &
 \includegraphics[width= 0.49 \linewidth]{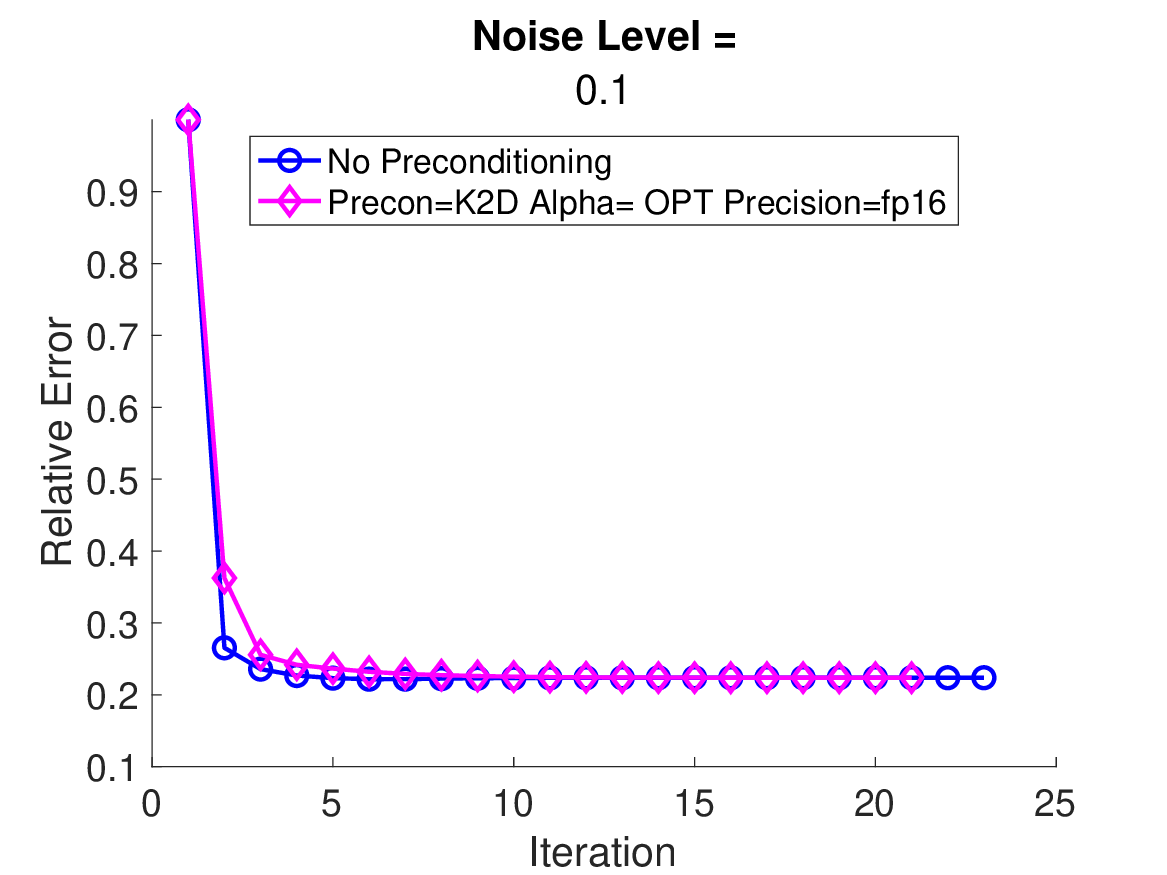} \\
 {(a)} & (b)
 \end{tabular}
 \end{center}
 \caption{A motion problem of size 128 (a) with 1\% noise; (b) with 10\% noise}
 \label{motion_128}
\end{figure}

\begin{figure}[!htb]
\begin{center}
\begin{tabular}{cc}
 \includegraphics[width= 0.49 \linewidth]{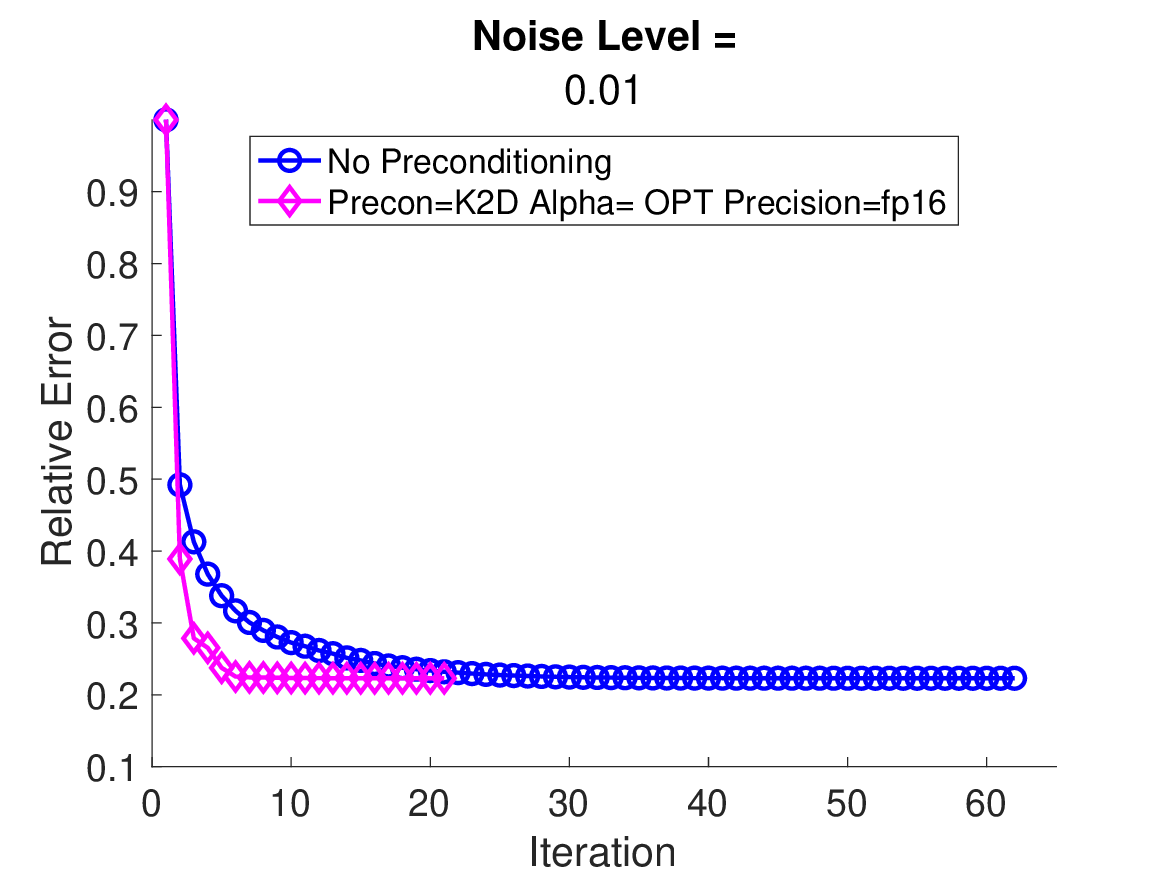} &
 \includegraphics[width= 0.49 \linewidth]{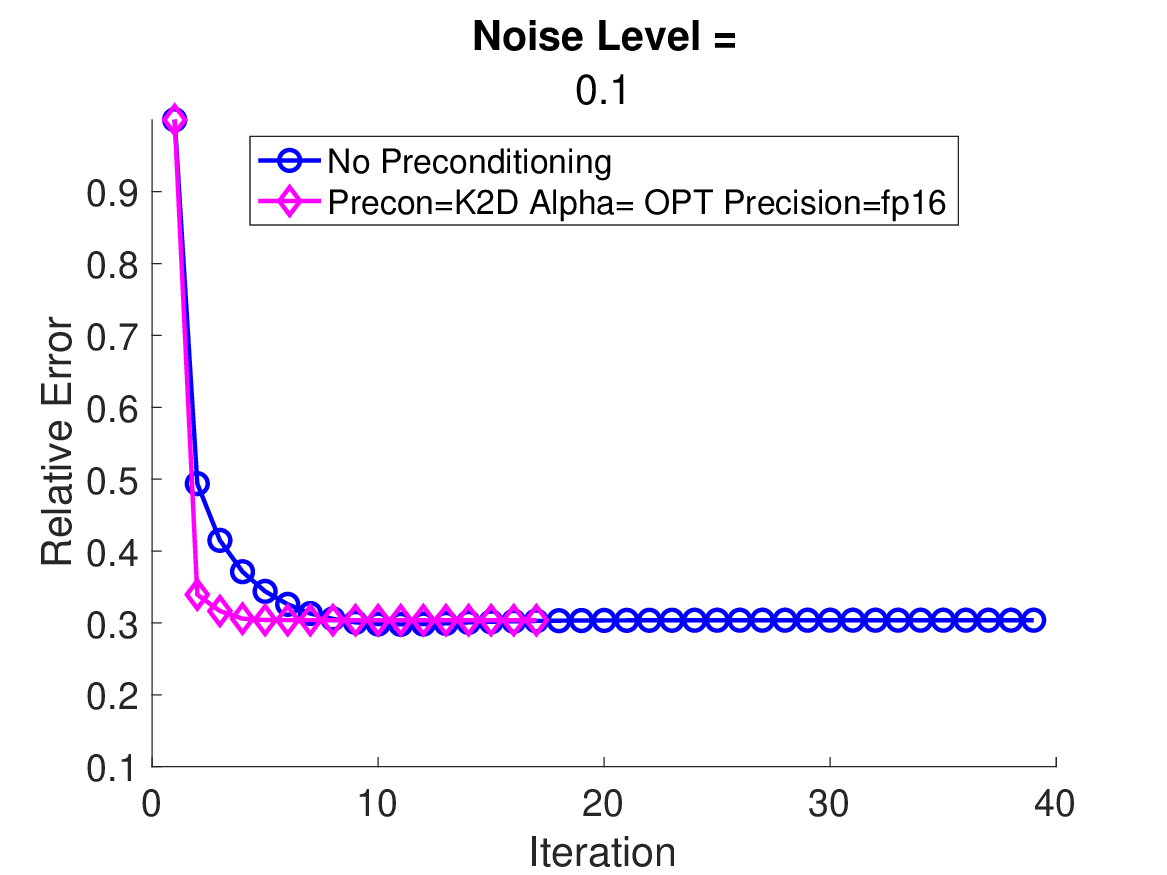} \\
 {(a)} & (b)
 \end{tabular}
 \end{center}
 \caption{A speckle problem of size 128 (a) with 1\% noise; (b) with 10\% noise}
 \label{speckle_128}
\end{figure}

The last blur we investigated is the speckle blur, which is used to model blur caused atmospheric turbulence \cite{gazzola2019irtools}. As summarized in Table~\ref{tab:PreconApprox}, in this case $A$ is a sum of $18$ Kronecker products. Similar to the shake blur, the preconditioning greatly improves the rate of convergence for low noise level, but the method is not very helpful when the noise level is high, as shown in Figure \ref{speckle_128}. Despite the large number of Kronecker product terms in $A$, our half-precision, 1-term approximation works effectively.

\section{Concluding Remarks}
\label{sec:ConcludingRemarks}

In this project, we investigated the use of a low-precision one-term Kronecker product approximation as a preconditioner for the PCG algorithm to solve ill-conditioned image deblurring problems. 
For the Gaussian blur, our preconditioning method converges significantly faster than the no preconditioning case since the Kronecker product approximation is exact, and thus the only approximation comes from using half-precision arithmetic. For other types of blur, our algorithm is more efficient when the noise level is low; with high noise, PCG has the same convergence rate as the no preconditioning method. Thus, our algorithm is more suitable for image blurring problems with low noise level. Additionally, in our experiments, we found that flexible PCG has similar performance as PCG. In any case, our low-precision algorithm is much cheaper than its double-precision counterpart without loss of accuracy, as discussed in Section \ref{workuniteval}, demonstrating the potentials of low-precision algorithms.

We also compared different regularization parameter choice methods, using our half-precision Kronecker product SVD approximation. While the normal GCV and discrepancy principle tend to choose $\lambda$ too small and blow up the relative error, we were able to get a reasonable parameter given a suitable weights ($\omega$ for wGCV and $\eta$ for the discrepancy principle). With suitably chosen weights, both methods can produce good regularization parameters, but we found that wGCV was less sensitive to the weight parameter. We manually chose an appropriate weight through trial and error in this project, and further research can be done on systematic ways of choosing a suitable regularization parameter.

\bibliographystyle{abbrv}
\bibliography{references}
\end{document}